\documentclass[11pt]{article}
\usepackage{amsmath,enumerate,fancyhdr,amssymb, amsthm, pstricks, epsf, lscape}
\usepackage{amsmath,enumerate,fancyhdr,amssymb, amsthm, lscape}
\usepackage{color,graphicx}
\usepackage{hyperref}
\usepackage[margin=1in]{geometry} 

\newcommand{\commentout}[1]{}
\newcommand{\co}[1]{}
\newcommand{\st}{\mbox{s.t. }}

\def\eref#1{(\ref{#1})}

\newcommand{\ti}{\times}

\newcommand{\norm}[1]{\parallel \! #1 \! \parallel}

\newcommand{\rad}[1]{\mathbb{R}^{#1}}

\newcommand{\zad}[1]{\mathbb{Z}^{#1}}

\newcommand{\latt}[1]{\mathbb{L}(#1)}
\newcommand{\nlatt}[1]{\mathbb{N}(#1)}

\newcommand{\LRA}{\Leftrightarrow}
\newcommand{\RA}{\Rightarrow}

\newcommand{\width}{{\rm width}}
\newcommand{\iwidth}{{\rm iwidth}}

\newcommand{\bb}{B\&B}

\newcommand{\la}{\langle}
\newcommand{\ra}{\rangle}

\newcommand{\nin}{\noindent}

\newcommand{\bx}{\bar{x}}

\newcommand{\tA}{\tilde{A}}
\newcommand{\ta}{\tilde{a}}

\newcommand{\tc}{\tilde{c}}
\newcommand{\tp}{\tilde{p}}
\newcommand{\tr}{\tilde{r}}
\newcommand{\hp}{\hat{p}}

\newcommand{\hr}{\hat{r}}

\newcommand{\tQ}{\tilde{Q}}

\newcommand{\hQ}{\hat{Q}}

\newcommand{\Frob}{\operatorname{Frob}}

\newcommand{\di}{\displaystyle}

\newtheorem{Definition}{Definition}

\newtheorem{Fact}{Fact}

\newtheorem{Setup}{Setup}
\newtheorem{Example}{Example}
\newtheorem{Counterexample}{Counterexample}
\newtheorem{Proposition}{Proposition}
\newtheorem{Lemma}{Lemma}
\newtheorem{Theorem}{Theorem}
\newtheorem{Corollary}[Definition]{Corollary}
\newtheorem{Remark}[Definition]{Remark}
\newtheorem{Assumption}{Assumption}
\newtheorem{Recipe}{Recipe}

\newcommand{\beq}{\begin{equation}}
\newcommand{\eeq}{\end{equation}}
\newcommand{\beqa}{\begin{eqnarray}}
\newcommand{\eeqa}{\end{eqnarray}}
\newcommand{\ba}{\begin{array}}
\newcommand{\ea}{\end{array}}
\newcommand{\bac}{\begin{array}{ccccccccccc}}
\newcommand{\eac}{\end{array}}
\newcommand{\bprop}{\begin{Proposition}}
\newcommand{\eprop}{\end{Proposition}}

\newcommand{\bcex}{\begin{Counterexample}}
\newcommand{\ecex}{\end{Counterexample}}

\newcommand{\beqast}{\begin{eqnarray*}}
\newcommand{\eeqast}{\end{eqnarray*}}
\newcommand{\benum}{\begin{enumerate}}
\newcommand{\eenum}{\end{enumerate}}
\newcommand{\bit}{\begin{itemize}}
\newcommand{\eit}{\end{itemize}}
\newcommand{\bth}{\begin{Theorem}}
\newcommand{\enth}{\end{Theorem}}

\newcommand{\bdef}{\begin{Definition}}
\newcommand{\Edef}{\end{Definition}}

\newcommand{\bsetup}{\begin{Setup}}
\newcommand{\esetup}{\end{Setup}}
\newcommand{\ble}{\begin{Lemma}}
\newcommand{\ele}{\end{Lemma}}
\newcommand{\bex}{\begin{Example}}
\newcommand{\eex}{\end{Example}}
\newcommand{\bcor}{\begin{Corollary}}
\newcommand{\ecor}{\end{Corollary}}
\newcommand{\brem}{\begin{Remark}}
\newcommand{\erem}{\end{Remark}}
\newcommand{\bass}{\begin{Assumption}}
\newcommand{\eass}{\end{Assumption}}

\newcommand{\brep}{\begin{Recipe}}
\newcommand{\erep}{\end{Recipe}}

\newcommand{\pf}[1]{\vspace{.35cm} \nin {\bf Proof {#1} }}

\newcommand{\bpx}{\begin{pmatrix}}
\newcommand{\epx}{\end{pmatrix}}

\newcommand{\bbx}{\begin{bmatrix}}
\newcommand{\ebx}{\end{bmatrix}}

\begin{document}

\title{\bf Column basis reduction \\ and decomposable knapsack problems}
\author{Bala Krishnamoorthy \thanks{Department of Mathematics, Washington State University, {\bf bkrishna@math.wsu.edu}} \\
G\'{a}bor Pataki   \thanks{Department of Statistics and Operations Research, UNC Chapel Hill, {\bf gabor@unc.edu}. Author supported by NSF award 0200308}
}
\date{}

\maketitle

\begin{abstract}
\noindent We propose a very simple preconditioning method for integer programming feasibility
problems: replacing the problem
$$
\ba{rcl}
b' \, \leq & Ax & \,\, \leq b \\
x & \in  & \zad{n}
\ea
$$
with
$$
\ba{rcl}
b' \, \leq & (AU)y & \,\, \leq b \\
y & \in  & \zad{n},
\ea
$$
where $U$ is a unimodular matrix computed via {\em basis reduction}, to make the columns of $AU$
short (i.e., have small Euclidean norm),
and nearly orthogonal (see e.g., \cite{LLL82}, \cite{K87}). 
Our approach is termed column basis reduction, and the reformulation is called 
rangespace reformulation. 
It is motivated by 
the technique proposed for
equality constrained IPs by Aardal, Hurkens and Lenstra. We also propose a simplified method 
to compute their reformulation.

We also study a family of IP instances, called {\em decomposable knapsack problems (DKPs)}. DKPs
generalize the instances 
proposed by Jeroslow, Chv\'atal and Todd, Avis, Aardal and Lenstra, and Cornu\'ejols et al. 
They are knapsack problems with a constraint vector of the form $pM + r, \,$ with $p >0$ and $r$ integral vectors, and $M$ a large integer.
If the parameters are suitably chosen in DKPs, we prove
\begin{itemize}
\item hardness results, when branch-and-bound branching on individual variables is applied;
\item that they are easy, if one branches on the constraint $px$ instead; and 
\item that branching on the last few variables in either the rangespace or the AHL reformulations is equivalent to 
branching on $px$ in the original problem. 
\eit
We also provide recipes to generate such instances.

Our computational study confirms that the behavior of the studied instances  in practice is as predicted by the theory.

\end{abstract}

\clearpage
\tableofcontents
\clearpage

\section{Introduction and  overview of the main results}
\label{sect-intro}

\paragraph[br]{Basis reduction} 
Basis reduction (BR for short) is a fundamental technique in 
computational number theory, cryptography, and integer programming.
If $A$ is a real matrix with $m$ rows, and $n$ independent columns, 
the {\em lattice} generated by the columns of $A$ is
\beq \label{def-latt-B}
\latt{A} \, = \, \{ \, Ax \, | \, x \in \zad{n} \, \}.
\eeq
The columns of $A$ are called a {\em basis} of $\latt{A}$. 
A square, integral matrix $U$ is {\em unimodular} if $\det U = \pm 1$.
Given $A$ as above, BR computes a unimodular $U$
such that the columns of $AU$ are ``short'' and ``nearly''
orthogonal. The following example illustrates the action of BR:
$$
A = \left( \!\! \begin{array}{cc} 289 & 18 \\
                   466 & 29 \\
                   273 & 17
    \end{array} \!\!  \right), \; U = \left( \!\! \begin{array}{cc}
                           1 & -15 \\ -16 & 241
                          \end{array} \!\!  \right),
\; AU = \left( \!\! \begin{array}{cc}   1 & 3 \\ 2 & -1 \\ 1 & 2
                                       \end{array} \!\!  \right).
$$
We have $\latt{A} = \latt{AU}. \,$ In fact for two matrices $A$ and $B$, 
$\latt{A} = \latt{B} \,$ holds, if and only if $B=AU \,$ for some $U$ unimodular matrix (see e.g. Corollary 4.3a, \cite{S86}).

In this work we use two BR methods. The first is the Lenstra, Lenstra, and Lov\'{a}sz (LLL for short) 
reduction algorithm  
\cite{LLL82} which runs in polynomial time for rational lattices. The second is
Korkhine-Zolotarev (KZ for short) reduction -- see \cite{K83} and \cite{Sh87} --
which runs in polynomial time for rational lattices only when the number of columns of $A$ is fixed. 

\paragraph[br-ip]{Basis reduction in Integer Programming} The first application of BR for integer programming is
in Lenstra's IP algorithm that runs in polynomial time in fixed dimension, see \cite{L83}. 
Later IP algorithms which share polynomiality for a fixed number of variables also 
relied on BR: see, for instance Kannan's algorithm \cite{K87}; Barvinok's algorithm to count the number of lattice 
points in fixed dimension \cite{Bar94, DK97}, and its variant proposed by de Loera et al.~in \cite{DHTY04}.
A related method in integer programming is {\em generalized basis reduction} due to Lov\'{a}sz and Scarf \cite{LS92}. For 
its implementation see Cook et. al in \cite{CRSS93}.  Mehrotra and Li in \cite{ML04} proposed a modification and implementation of Lenstra's method, and of 
generalized basis reduction. For surveys, we refer to \cite{S86} and \cite{K87agn}. 

A computationally powerful reformulation technique based on BR 
was proposed by Aardal, Hurkens, and Lenstra in \cite{AHL00}. They 
reformulate an equality constrained IP feasibility problem
\beq \label{eq-ip}
\ba{rcl}
Ax & = & b \\
\ell \leq & x & \leq u \\
       & x & \in \zad{n}
\ea
\eeq
with integral data, and $A$ having $m$ independent rows,  as follows: they find a matrix
$B$, and a vector $x_b$ with $[B, x_b]$ having short, and nearly
orthogonal columns, $x_b$ satisfying $\, A x_b = b, $ and the property
\beq \label{aardal-ref}
\{ \, x \in \zad{n} \, | \, A x = 0 \, \} = \{ \, B \lambda  \, | \, \lambda \in \zad{n-m} \, \}.
\eeq

The reformulated instance is 
\beq \label{eq-ip-ref}
\ba{rcl}
\ell - x_b \leq & B \lambda   & \leq u - x_b \\
       & \lambda & \in \zad{n-m}.
\ea
\eeq
For several families of hard IPs, the reformulation \eref{eq-ip-ref} 
turned out to be much easier 
to solve for commercial MIP solvers than the original one; a notable family 
was the {\em marketshare} problems of Cornu\'ejols and Dawande \cite{CD98}. The solution 
of these instances using the above reformulation technique is described 
by Aardal, Bixby, Hurkens, Lenstra, and Smeltink in \cite{ABHLS00}.

The matrix $B$ and the vector $x_b$ are found as follows.
Assume that $A$ has $m$ independent rows. They embed $A$
and $b$ in a matrix, say $D$, with $n+m+1$ rows, and $n+1$ columns,
with some entries depending on two large constants $N_1$ and
$N_2$: 
\beq \label{def-D}
D = \bpx I_{n} & 0_{n \ti 1} \\
         0_{1 \ti n}   & N_1 \\
         N_2 A & - N_2 b
     \epx.
\eeq
The lattice generated by $D$  looks like
\beq \label{lattD} 
\latt{D} \, = \, \left\{ \left( \ba{c} x \\ N_1 x_0 \\ N_2(Ax-b x_0) \ea \right) \, \left| \right. \, \left( \ba{c} x \\ x_0  \ea \right) \in \zad{n+1} \, \right\},
\eeq
in particular, all vectors in a reduced basis of $\latt{D}$ have this form. 

For instance, if $A = [ 2, 2, 2], \,$ and $ b = 3, \,$ (this corresponds to the infeasible IP $2 x_1 + 2 x_2 + 2 x_3 = 3, \, x_i \in \{0,1\}, \,$ when the bounds
on $x$ are $0$ and $e$), then 
\beq \label{lattD2} 
\latt{D} \, = \, \left\{ \left( \ba{c} x \\ N_1 x_0 \\ N_2( 2 x_1 + 2 x_2 + 2 x_3 - 3 x_0) \ea \right) \, \left| \right. \, \left( \ba{c} x \\ x_0  \ea \right) \in \zad{3} \times \zad{1} \, \right\}.
\eeq
It is shown in \cite{AHL00},
that if $N_2 >\!\!> N_1 >\!\!> 1$ are suitably chosen, then in a reduced basis of $\latt{D}$ 
\bit
\item  $n - m \,$ vectors arise from some $\bpx x \\ x_0 \epx$ with
$ 
Ax= b x_0, \, x_0 = 0,
$ and 
\item $1 \,$ vector will arise from an $\bpx x \\ x_0 \epx$ with
$
Ax=b x_0, \, x_0 = 1.
$
\eit
So the $x$ vectors from the first group can form the columns of $B$, and the $x$ from the last can 
serve as $x_b$. If LLL- or KZ-reduction (the precise definition is given later) is used to compute the 
reduced basis of $\latt{D}, \,$ then $B$ is a basis reduced in the same sense. 

Followup papers on this reformulation technique were written by 
Louveaux and Wolsey \cite{LW02}, and Aardal and Lenstra \cite{AL04, AL06}. 

\paragraph[ques]{Questions to address} The speedups obtained by the Aardal-Hurkens-Lenstra (AHL) reformulation lead to 
the following questions:
\benum
\item[(Q1)] Is there a similarly effective reformulation technique for 
{\em general} (not equality constrained) IPs? 
\item[(Q2)]  Why does the reformulation work? Can we analyse its action on a reasonably wide class
of difficult IPs?
\eenum
More generally, one can ask:
\benum
\item[(Q3)]  What kind of integer programs are hard for a certain standard approach, such as branch-and-bound branching on individual variables,
and easily solvable by a different approach? 
\eenum
As to (Q1), one could simply add slacks to turn inequalities into equalities, and then apply the AHL reformulation.
This option, however, has not been studied. The
mentioned papers emphasize the importance of reducing the dimension of
the space, and of the full-dimensionality of the reformulation. Moreover, 
reformulating an IP with $n$ variables, $m$ dense constraints, and some bounds in this way
leads to a $D$ matrix (see \eref{def-D}) with $n+2m+1$ rows and $n+m+1$ columns.

A recent paper of Aardal and Lenstra \cite{AL04}, and  \cite{AL06} 
addressed the second question. They considered an 
equality constrained knapsack problem with unbounded variables
\beq  \label{kp2-eq} \tag\mbox{{KP-EQ}}
  \ba{rl}
ax & \,= \beta \\
x  & \, \geq 0 \\
x & \, \in \zad{n},
  \ea
  \eeq
with the constraint vector $a$ decomposing as $a = pM+r, \, $ with $p, r \in \zad{n}, \,$ $p>0, \,$ $M$ a positive integer, 
under the following  assumption: 
\bass \label{al-assumption} 
\benum
\item $r_j/p_j = \max_{i = 1, \dots, n } \, \{ r_i/p_i \}, \, r_k/p_k = \min_{i = 1, \dots, n} \, \{ r_i/p_i \}.$
\item $a_1 < a_2 < \dots < a_n$;
\item $\sum_{i=1}^n |r_i| < 2M$;
\item $M > 2 - r_j/p_j$; 
\item $M > r_j/p_j - 2r_k/p_k$.
\eenum
\eass
They proved the following:
\benum
\item Let $\Frob(a)$ denote the Frobenius number of $a_1, \dots, a_n, \,$ 
i.e., the  largest $\beta$ integer for which 
\eref{kp2-eq} is infeasible. Then 
\beq \label{AL-lb}
\Frob(a) \geq  \dfrac{(M^2p_j p_k +M(p_j r_k  + p_k r_j) +r_j r_k)(1 - \dfrac{2}{M+r_j/p_j})}{p_k r_j - p_j r_k} - (M + r_j/p_j).
\eeq
\item In the reformulation \eref{eq-ip-ref}, if we denote the last column of $B$ by $b_{n-1}, \,$  then 
\beq \label{AL-lb2}
\norm{b_{n-1}} \geq \frac{\norm{a}}{\sqrt{ \norm{p}^2\norm{r}^2 - (pr^T)^2 }}.
\eeq
\eenum
It is argued in \cite{AL04} that 
the large right-hand side explains the hardness of the corresponding instance, and 
that the large norm of $b_{n-1}$ explains why the reformulation is easy: if we branch on $b_{n-1}$ in the 
reformulation, only a small number of nodes are created in the branch-and-bound tree. 
In  section \ref{critique} we show that there is a gap in the proof of \eref{AL-lb2}. Here we also show an instance 
of a bounded polyhedron where the columns of the  constraint matrix are LLL-reduced, but branching on a variable corresponding to the longest column
produces exponentially many nodes. 

Among the other papers that motivated this research, two are ``classical'': \cite{J74}, and parts of \cite{C80}. 
They all address the hardness question in (Q3), and the easiness is straightforward to show. 

Jeroslow's knapsack instance in \cite{J74} is 
\beq \label{jeroslow}
\ba{rl}
\min & x_{n+1}  \\
st.  & 2 \sum_{i=1}^n x_i + x_{n+1} \, = \, n   \\
     & x_i \in  \{ 0, 1 \} \, (i=1, \dots, n+1),
\ea
\eeq
where $n$ is an odd integer. The optimal solution of \eref{jeroslow} is trivially $1$, but 
branch-and-bound requires an exponential number of nodes to prove this, if we branch on individual variables. 

In \cite{C80} Todd and Avis constructed knapsack problems
of the form 
\beq
\ba{rl} 
\max & ax \\
st. & ax \leq \beta \\
     & x \in \{ 0, 1 \}^n,
\ea
\eeq
with $a$ decomposing as $a = eM+r  \,$ ($M$ and $r$ are chosen differently in the Todd- and in the Avis-problems). 
They showed that these instances exhibit a similar behavior. 
\co{Though the latter are optimization problems, they are hard for the same reason as their feasibility version with 
$ax = \beta$ is: even with many $x_i$ variables fixed, there is a solution which is feasible to the LP-relaxation. }

Though this is not mentioned in \cite{J74}, or \cite{C80}, it is straightforward to see that the Jeroslow-, Todd-, and Avis-problems 
can be solved at the rootnode, if one branches on the constraint $\sum_{i=1}^n x_i \,$ instead of branching on the $x_i$. 

A more recent work that motivated us is \cite{CUWW97}. Here a family of instances 
of the form 
\beq \label{cu}
\ba{rl} 
\max & ax \\
\st & ax \leq \beta \\
     & x \in \zad{n}_+,
\ea
\eeq
with $a$ decomposing as $a = pM+r, \, $ where $p$ and $r$ are integral vectors, and $M$ is a positive integer, 
was proposed.  The authors used $\Frob(a)$ as $\beta \,$ in \eref{cu}. These problems turned out to be hard for commercial MIP solvers, 
but easy if one  uses a test-set approach. One can also verify computationally that if one branches on the constraint $px$ in these instances, 
then feeds the resulting subproblems to a commercial solver, they are solved quite quickly.

\paragraph[contrib]{Contributions, and organization of the paper} 

We first fix basic terminology. When branch-and-bound (\bb \ for short) branches on individual variables, we call the resulting algorithm 
{\em ordinary branch-and-bound}. 

\bdef \label{basicdef}
If $p$ is an integral vector, and $k$ an integer, then the logical expression 
$px \leq k \, \vee \, px \geq k+1$ is called a {\em split disjunction}. 
We say that the infeasibility of an integer programming problem is proven 
by $px \leq k \, \vee \, px \geq k+1, \,$ if both polyhedra 
$\{ \, x \, | \, px \leq k \, \}$ and $\{ \, x \, | \, px \geq k+1 \, \}$ 
have empty intersection
with the feasible set of its LP relaxation.

We say that the infeasibility of an integer programming problem is proven by branching on $px, \,$ if
$px \,$ is nonintegral for all $x$ in its LP relaxation.

\end{Definition} 

We call  a knapsack problem with weight vector $a$ a {\em decomposable knapsack problem} (DKP for short), 
if $a = pM+r, \,$ where $p \,$ and $r$ are integral vectors, $p>0, \,$ and $M $ is a large integer. 
We could not find a good definition of DKPs which would not be either too restrictive, or too permissive, as far as how large $M$ should be.
However, we will show how to find $M$ and the bounds for given $p$ and $r$ so the resulting DKP has interesting properties.

\co{for generate DKP instances with interesting properties. 
We do not define DKPs formall
, and 
$\norm{a} > \norm{r}$ and $M > \norm{r}$ hold.  
}

The paper focuses on the interplay of these concepts, and their connection to IP reformulation techniques. 

\benum

\item In the rest of this section we describe a simple reformulation technique, called the {\em rangespace reformulation} for 
arbitrary integer programs. The dimension of the reformulated 
instance is the same as of the original. We also show a simplified method to compute the AHL reformulation, and 
illustrate how the reformulations work on some simple instances. 

For a convenient overview of the paper we state Theorems \ref{summ1} and \ref{summ2} as a sample of the main results. 

\item In Section \ref{dkp-hard} we consider knapsack feasibility problems with a positive weight vector. We show a somewhat surprising result: 
if the infeasibility of such a problem is proven by $px \leq k \vee px \geq k+1, \,$ with $p$ positive,
then  a lower bound follows on the number of nodes that must be enumerated by ordinary \bb \ to prove infeasibility.
So, easiness for constraint branching implies hardness for ordinary \bb .

\item In Section \ref{recipes} we give two recipes to find  DKPs, whose infeasibility is proven by the split disjunction $px \leq k \vee px \geq k+1. \,$ 
Split disjunctions for deriving cutting planes have been studied e.g. in \cite{CKS90, CL03, B79, BCC93}. This paper seems to be 
the first systematic study of knapsack problems with their infeasibility having such a short certificate. 

Thus (depending on the parameters), their hardness for ordinary \bb \ follows using the results of Section \ref{dkp-hard}. 
We show that several well-known hard integer programs from the literature, such as the
Jeroslow-problem \cite{J74}, and the Todd- and Avis-problems from \cite{C80} can be found using Recipe 1. 
Recipe 2 generates instances of type \eref{kp2-eq}, with a short proof (a split disjunction) of their infeasibility. 

So this section provides a unifying framework to show 
the hardness of instances (for ordinary \bb ) which are easy for constraint branching. These results add to the understanding of hard knapsacks 
described in \cite{J74}, \cite{C80} and \cite{AL04}, as follows.
We deal with arbitrary knapsacks, both with bounded, 
and unbounded variables; we give explicit lower bounds on the number of nodes that ordinary \bb \ must enumerate, which is done in \cite{C80} 
for the Todd and Avis instances; and our instances have a short, split disjunction certificate. 

Using the recipes we generate some new, interesting examples. For instance,
Example \ref{kp2-eq-nt} is a knapsack problem whose infeasibility is proven by a single split disjunction, but 
ordinary \bb \ needing a {\em superexponential} number of nodes to prove the same. Example \ref{reverse-avis} reverses the role of the two vectors 
in the Avis-problem, and gives an instance which is computationally more difficult than the original.

\item In Section \ref{sect-large-rhs} we extend the lower bound \eref{AL-lb} in two directions. 
We first show that for given $p $ and $r$ integral vectors,
and  sufficiently large $M, \,$ 
there is a range of $\beta$ integers for which the infeasibility of 
\eref{kp2-eq} with $a = pM+r \,$ is proven by branching on $px$. The smallest such integer is essentially the same 
as the lower bound in  \eref{AL-lb}. 

Any such $\beta$ right-hand side is a lower bound on $\Frob(a), \,$ 
with a short certificate of being a lower bound, i.e. a split disjunction certificate of the infeasibility of \eref{kp2-eq}.

We then study the {\em largest} integer for which the infeasibility of \eref{kp2-eq} with $a = pM+r,  \,$ and $M$ sufficiently large, 
is proven by branching on $px$. We call this number the $p$-branching Frobenius number, and give 
a lower and an upper bound on it. 

\item In Section \ref{geometry} we show some basic results on the geometry of the reformulations. 
Namely, given a vector say  $c, \,$ we find a vector which achieves the same width in the reformulation, as $c$ does in the original problem.

\item Subsection \ref{subsec-range}  shows why DKPs 
become easy after the rangespace reformulation is applied. 
In Theorem \ref{range-main-thm} 
 we prove that if $M$ is sufficiently large, and the infeasibility of a DKP is proven by branching on $px, \,$ then 
the infeasibility of the reformulated problem is proven by  branching on the last few variables in the reformulation. 
How many ``few'' is will depend on the magnitude of  $M$. 
We give a similar analysis for the AHL reformulation in Subsection \ref{subsec-null}. \co{, and here we also point out the gap in 
the proof of the lower bound \eref{AL-lb2}. }

Here we remark that a method which
{\em explicitly} extracts ``dominant'' directions in an integer program was proposed by
Cornu\'ejols at al in in \cite{CUWW97}.

\item In Section \ref{comp} we present a computational study that compares 
the performance of an MIP solver before and after the application of the reformulations on certain 
DKP classes. 

\item In Section \ref{critique} we point out a gap in the proof of \eref{AL-lb2}, and show a correction. 
We also describe a bounded polyhedron with the columns of the  constraint matrix forming an LLL-reduced basis, 
where branching on a variable corresponding to the longest column creates exponentially many subproblems.
\eenum

\paragraph[ref]{The rangespace reformulation} 

Given 
\beq \label{ip} \tag{IP}
\ba{rcl}
b' \, \leq & Ax & \,\, \leq b \\
x & \in  & \zad{n},
\ea
\eeq
we compute a unimodular (i.e. integral, with $\pm 1$ determinant) matrix $U$ that makes
the columns of $AU$  short, and nearly orthogonal; 
$U$ is computed using basis reduction, either the LLL- or the KZ-variant 
(our analysis will be unified).
We then recast \eref{ip} as
\beq \label{tip} \tag{$\tilde{\mbox{IP}}$}
\ba{rcl}
b' \, \leq & (AU) y & \,\, \leq b \\
y & \in  & \zad{n}.
\ea
\eeq
The dimension of the problem is unchanged; we will call this technique
{\em rangespace reformulation}. 

\bex \label{ex1}
{\rm
Consider the infeasible problem
\beq \label{ex1-prob}
\ba{rcl}
106 \, \leq & 21 x_1 + 19 x_2 & \leq \, 113 \\
0   \, \leq & x_1, x_2 & \leq \, 6 \\
            & x_1, x_2 \in \zad{}, &
\ea
\eeq
with the feasible set of the LP-relaxation depicted on the first picture in Figure \ref{fig-ex1}.
In a sense it is both hard, and easy. 
On the one hand, branching on either variable will produce at least $5$ feasible nodes. 
On the other hand, the maximum and the minimum of $x_1 + x_2$ over the LP relaxation of
\eref{ex1-prob} are $5.94$, and $5.04$, respectively, thus ``branching''
on this constraint proves infeasibility at the root node.

\begin{figure}
\label{fig-ex1}
\epsfxsize=3in
\epsffile{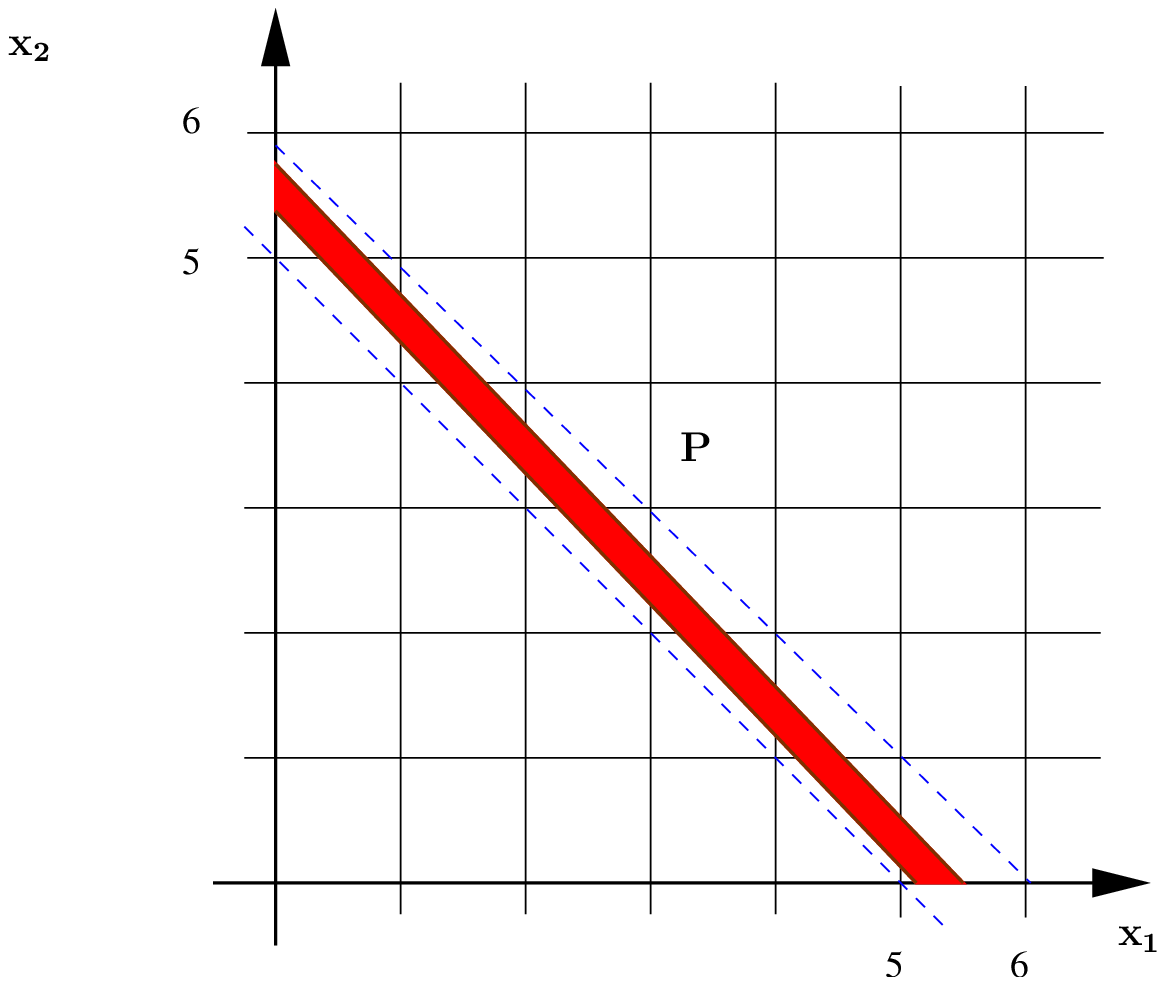}
\hspace{0.01in}
\epsfxsize=3in
\epsffile{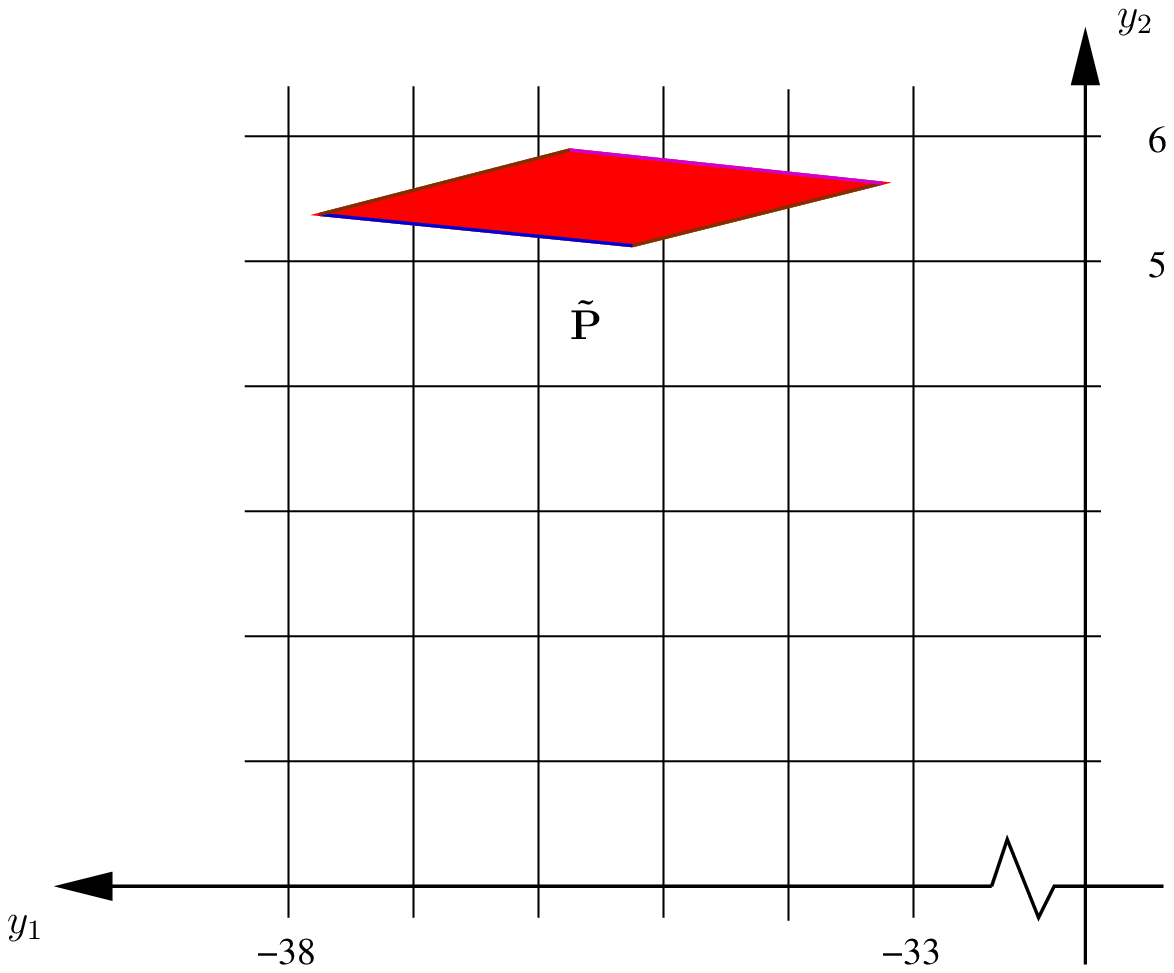}
\caption{The polyhedron in Example \ref{ex1} before and after reformulation}
\end{figure}
When the rangespace reformulation with LLL-reduction is applied, we have
$$
A = \bpx 21 & 19 \\ 1 & 0 \\ 0 & 1 \epx, \, U = \bpx -1 & - 6 \\ 1 & 7 \epx, \, AU = \bpx -2 & 7 \\ -1 & - 6 \\ 1 & 7 \epx,
$$
so the reformulation is
\beq \label{ex1-ref}
\ba{rcl}
106 \, \leq & -2 y_1 + 7 y_2 & \leq \, 113 \\
0   \, \leq & - y_1 - 6 y_2  & \leq \, 6 \\
0   \, \leq & y_1 + 7 y_2  & \leq \, 6 \\
            & y_1, y_2 \in \zad{}.   \\
\ea
\eeq
Branching on $y_2$ immediately proves infeasibility, as the second picture in Figure \ref{fig-ex1} shows.
The linear constraints of \eref{ex1-ref} imply
\beq \label{y2-504}
5.04 \leq y_2 \leq 5.94.
\eeq
These bounds are the same as the bounds on $x_1 + x_2$.  This fact will follow 
from Theorem \ref{range-geometry}, a general result about how the widths 
are related along certain directions in the original and the reformulated problems. 
}
\eex

\bex \label{ex2} {\rm This example is a simplification of Jeroslow's problem \eref{jeroslow} 
from \cite{J74}. Let $n$ be a positive odd integer. 
The problem
\beq \label{ex2-prob}
\ba{rcl}
 2 \sum_{i=1}^n x_i & = & n \\
 0 \leq  x   & \leq & e \\
    x & \in & \zad{n} 
\ea
\eeq
is integer infeasible. 

Ordinary \bb \ (i.e. \bb \ branching on the $x_i$ variables) must enumerate at least $2^{(n-1)/2}$ nodes
to prove infeasibility. To see this, suppose that at most $(n-1)/2$ variables are fixed to either $0$ or $1$.
The sum of the coefficients of these variables is at most $n-1$, while the sum of the 
coefficients of the free variables
is at least $n+1$. Thus, we can set some free variable(s) to a possibly 
fractional value to get an LP-feasible solution.

On the other hand, denoting by $e$ the vector of all ones, 
the maximum and minimum of $ex$ over the LP relaxation of \eref{ex2-prob} is
$n/2, $ thus branching on $ex$ proves infeasibility at the root node.

For the rangespace reformulation using LLL-reduction, we have
$$
A = \bpx 2e_{1 \ti n} \\ I_n \epx, \, U = \bpx I_{n-1} & 0_{(n-1) \ti 1} \\ -e_{1 \ti (n-1)} & 1 \epx, \, AU = \bpx 0_{1 \ti (n-1)} & 2 \\  I_{n-1} & 0_{(n-1) \ti 1} \\ -e_{1 \ti (n-1)} & 1 \epx, \
$$
thus the reformulation is
\beq \label{ex2-ref}
\ba{rcl}
            &  2 y_n & = n  \\
     0 \leq & y_1, \dots, y_{n-1}  & \leq 1 \\
     0 \leq & - \sum_{i=1}^{n-1} y_i + y_n & \leq 1 \\
                     & y \in \zad{n}. & \\
\ea
\eeq
So branching on $y_n$ immediately implies the infeasibility of \eref{ex2-ref}, and thus of
\eref{ex2-prob}. 
}
\eex

\paragraph[col-br-null]{A simplified method to compute the AHL reformulation} 
Rangespace reformulation only affects the constraint matrix,
so it can be applied unchanged, if some of the two-sided inequalities in \eref{ip} are
actually equalities, as in Example \ref{ex2}. We can still choose a different way of reformulating the problem. 
Suppose that 
\beq \label{subsys}
A_1 x = b_1
\eeq
is a system of equalities contained in the constraints of \eref{ip}, and assume that
$A_1$ has $m_1$ rows.
First compute an integral matrix $B_1$, 
$$
\{ \, x \in \zad{n} \, | \, A_1 x = 0 \, \} = \{ \, B_1 \lambda \, | \, \lambda \in \zad{n-m_1} \, \},
$$
and an integral vector $x_1$ with $A x_1 = b_1$. 
$B_1$ and $x_1$ can be found by a Hermite Normal Form computation -- see e.g. \cite{S86}, page 48.

In general, the columns of $[B_1, x_1]$ will {\em not} be reduced.
So, we substitute $B_1 \lambda + x_1$ into the part of \eref{ip} excluding \eref{subsys}, and
apply the rangespace reformulation to the resulting system.

If the system \eref{subsys} contains all the constraints of the integer program other than the bounds, then 
this way we get the AHL reformulation. 

\bex \label{ex3}
{\rm {\bf (Example \ref{ex2} continued)} In this example 
\eref{ex2-prob} has no solution over the integers, irrespective of the bounds. 

However,  we can rewrite it as 
\beq \label{ex2-prob-eq}
\ba{rcl}
2 \sum_{i=1}^n x_i + x_{n+1} & = & n  \\
0 \,\,\, \leq  x_{1:n}   & \leq & e  \\
- 1/2 \,\,\, \leq x_{n+1} & \leq & 1/2 \\
 x \in \zad{n}. & \\
\ea
\eeq
The $x$ integer vectors that satisfy the first equation in \eref{ex2-prob-eq} can be parametrized  with $\lambda \in \zad{n} \,$ as 
\beq \label{xlambda}
\ba{rcl}
x_1 & = & \lambda_1 + \dots + \lambda_n \\
x_2 & = & - \lambda_1 \\
    & \vdots &        \\
x_n &  = & - \lambda_{n-1} \\
x_{n+1} & = & - 2 \lambda_n + n.
\ea
\eeq
Substituting \eref{xlambda} into the bounds of \eref{ex2-prob-eq} we obtain the reformulation 
\beq \label{ex2-prob-eq-ref}
\ba{rcl}
0 \leq & \sum_{j=1}^{n-1} \lambda_j + \lambda_n & \leq 1 \\
0 \leq & - \lambda_j                            & \leq 1 \; (j=1, \dots, n-1) \\
- 1/2 \leq & -2 \lambda_n +n                    & \leq 1/2 \\
       & \lambda \in \zad{n}. & 
\ea
\eeq
The columns of the constraint matrix of \eref{ex2-prob-eq-ref} are already reduced in the LLL-sense. 
The last constraint is equivalent to 
\beq \label{lambdaninf}
(n+1)/2 - 3/4 \leq \lambda_n \leq (n+1)/2 - 1/4,
\eeq
so the infeasibility of \eref{ex2-prob-eq-ref} and thus of \eref{ex2-prob-eq}
is proven by branching on $\lambda_n$.
}
\eex

\paragraph[cbr-rhs]{Right-hand side reduction}
On several instances we found that reducing the right-hand side in \eref{ip} yields
an even better reformulation.  To do this, we rewrite \eref{ip} as
\beq \label{ip2} \tag{IP2}
\ba{rl}
 F x & \, \leq f \\
x & \in \, \zad{n},
\ea
\eeq
then reformulate the latter as
\beq \label{tip2} \tag{$\tilde{\mbox{IP2}}$}
\ba{rl}
(FU)y & \, \leq f  - (FU) x_r \\
y     & \, \in  \zad{n},
\ea
\eeq
where the unimodular $U$ is again computed by basis reduction,
and $x_r \in \zad{n}$ to make $f - (FU) x_r$ short, and near orthogonal to the columns of
$FU$. For the latter task, we may use --  for instance --
 Babai's algorithm \cite{Ba86} to find $x_r, \, $ so that $(FU) x_r$ is a nearly closest vector to $f$
in the lattice generated by the columns of $F$.

It is worth to do this, if the original constraint matrix, and right-hand side (rhs) {\em both} have 
large numbers. Since the rangespace reformulation reduces the matrix coefficients, leaving large numbers in the rhs may lead 
to numerical instability. Our analysis, however, will rely only on the reduction of the constraint matrix.

\paragraph[ahl-comp]{Rangespace, and AHL reformulation}
To discuss the connection of these techniques, we assume for simplicity that right-hand-side reduction 
is not applied. 

Suppose that $A \,$ is an integral matrix with 
$m$ independent rows,  and $b \,$ is an integral column vector with $m$ components. Then 
the equality constrained IP 
\beq \label{eq-ip-2}
\ba{rcl}
Ax & = & b \\
\ell \leq & x & \leq u \\
       & x & \in \zad{n}
\ea
\eeq
has another, natural formulation:
\beq \label{eq-ip-ref-2}
\ba{rcl}
\ell \leq & B \lambda + x_b  & \leq u \\
       & \lambda & \in \zad{n-m},
\ea
\eeq
where 
\beq \label{ref-2}
\{ \, x \in \zad{n} \, | \, A x = 0 \, \} = \{ \, B \lambda  \, | \, \lambda \in \zad{n-m} \, \},
\eeq
and $x_b \,$ satisfies $A x_b = b. \,$ The matrix $B$ can be constructed from $A$ using an HNF computation. 

Clearly, to \eref{eq-ip-2} we can apply 
\bit
\item the rangespace reformulation (whether the constraints are inequalities, or equalities), or 
\item  the AHL method, which is equivalent to applying the rangespace reformulation to \eref{eq-ip-ref-2}. 
\eit
So, on  \eref{eq-ip-2} the rangespace reformulation method can be viewed as
a ``primal'' and the AHL reformulation as a ``dual'' method. The somewhat surprising fact is, 
that for a fairly large class of problems {\em both}  work, both theoretically, and computationally.
When both methods are applicable, 
we did not find a significant difference in their performance on the tested problem instances.

An advantage of the rangespace reformulation 
is its simplicity. For instance, there is a one-to-one correspondence between ``thin''
branching directions in the original, and the reformulated problems, so 
in this sense the geometry of the feasible set is preserved.  The correspondence is described in Theorem \ref{range-geometry} 
in Section \ref{geometry}. The situation is more complicated for the AHL method, and 
correspondence results are described in Theorems \ref{null-geometry} and \ref{null-geometry-reverse}.
These  results use ideas from, and generalize Theorem 4.1 in \cite{ML04}. 

In a sense the AHL method can be used to simulate the rangespace method on an 
{\em inequality constrained } problem: we can simply add slacks beforehand. 
However: 
\bit
\item the rangespace reformulation can be applied to an equality constrained problem as well, 
where there are no slacks; 
\item the main point of our paper is not simply presenting a 
reformulation technique, but {\em analysing} it. The analysis must be carried out separately for 
the rangespace and AHL reformulations. In particular, the bounds on $M$ that ensure 
that branching on the ``backbone'' constraint  $px$ in \eref{kp2} will be 
mimicked by branching on a small number of individual variables in the reformulation will be smaller
in the case of rangespace reformulation.
\eit

Using the rangespace reformulation is also natural when dealing with an optimization problem
of the form 
\beq \label{ipopt} \tag{IP-OPT}
\ba{rrrcl}
\max & & c x \\
\st  & b' \leq & Ax & \leq & b \\
     &  & x & \in  & {\mathbb Z}^n.
\ea
\eeq
Of course, we can reduce solving \eref{ipopt} to a sequence of feasibility problems. 

A simpler method is solving \eref{ipopt} by {\em direct reformulation}, i.e. by solving 
\beq \label{refipopt} \tag{$\widetilde{\mbox{IP-OPT}}$} 
\ba{rrrcl}
\max & & \tc  y \\
st.  & b' \leq & \tA y & \leq & b \\
     &         & y & \in  & {\mathbb Z}^n,
\ea
\eeq
where 
$$
\tc = cU, \, \tA = AU,
$$ 
with $U$ having been computed to make the columns of 
$$
\bpx  c  \\
      A  
\epx U 
$$
reduced. 

\paragraph[ref]{Some other reformulation methods} 
Among early references, the all-integral simplex algorithm of Gomory \cite{G60} can  be viewed as a reformulation method.
Bradley in \cite{Brad71} studied integer programs connected via 
unimodular transformations, akin to how the rangespace reformulation works. 
However, the transformations in \cite{Brad71} do not arise from basis reduction.

The Integral Basis Method \cite{HKW03} has two reformulation steps: in the first an integral basis 
of an IP from a nonintegral basis of the LP relaxation is found. In the second, an augmentation vector leading to a better integral solution is found, or shown not to exist. 
Haus in his dissertation \cite{Hau04} studied the question of how to derive such augmentation vectors for general IPs. 

\paragraph[not]{Notation} Vectors are denoted by lower case letters. In notation we do not distinguish 
between row and column vectors; the distinction will be clear from the context. Occasionally, we write 
$\la x, y \ra$ for the inner product of vectors $x$ and $y$. 

We denote the sets of nonnegative, and positive integers by 
$\zad{}_+, \, $ and $\zad{}_{++},  $ respectively. The sets of nonnegative, and positive integral $n$-vectors are denoted 
by $\zad{n}_+, \, $ and $\zad{n}_{++},$ respectively. If $n$ a positive integer, then  $N$ is the set $\{ 1, \dots, n \}$. 
If $S$ is a subset of $N, \,$ and $v$ an $n$-vector, then $v(S)$ is defined as $\sum_{i \in S} v_i$. 

For a matrix $A$ we use a Matlab-like notation, 
and denote its $j$th row, and column by 
$
A_{j, :} \; \text{and} \; A_{:, j},
$
respectively. Also, we denote the subvector $(a_k, \dots, a_\ell)$ of a vector
$a$ by $a_{k : \ell}$. 

For $p \in \zad{n}_{++}, \,$ and an integer $k$ we write 
\beq \label{defratio} 
\ba{rcl}
\ell(p, k) & = & \max \, \{ \, \ell \, | \, p(F) \leq k, \, \text{and} \, p(N \setminus F) \geq k+1 \, \forall \, F \subseteq N, \, |F| = \ell \, \}.
\ea
\eeq
The definition implies that $\ell(p,k)=0 \,$ if $k \leq 0, \,$ or $k \geq \sum_i p_i, \,$ and 
$\ell(p,k)$ is large if the components of $p$ are small relative to $k, \,$ and not too different from each other. 
For example, if $p = e, \, k < n/2, \,$ then $\ell(p, k) = k.$ 

Sometimes $\ell(p,k)$ is not easy to compute exactly, but we can use a good lower bound, which is usually easy to find.
For instance, let $n$ be an integer divisible by $4, \,$ $p = (1, 2, \dots, n). \,$ 
The first $3n/4$ components of $p$ sum to strictly more than 
$(\sum_{i=1}^n p_i)/2,$ and the last $n/4$ sum to strictly less than this. 
Since the components of $p$ are ordered increasingly, it follows that 
$$
\ell(p, n(n+1)/4) \, \geq \, n/4.
$$
On the other hand, $\ell(p,k)$ can be zero, even if $k$ is positive. For example,
if $p \,$ is {\em superincreasing}, i.e. $p_i > p_1 + \dots + p_{i-1} \,$ for $i = 2, \dots, n, \,$ then it is easy to see that 
$\ell(p, k) = 0 \,$ for any positive integer $k$. 

\paragraph[dkp]{Knapsack problems} We will study knapsack feasibility problems 
 \beq  \label{kp2} \tag{KP}
  \ba{rcl}
\beta_1 \, \leq & ax & \,\, \leq \beta_2 \\
0 \, \leq & x & \,\, \leq u \\
   &  x & \, \in \zad{n}.
  \ea
  \eeq
In the rest of the paper for the data of \eref{kp2} we will use the following assumptions, that we collect here for convenience: 
\bass
The row vectors $a, u \,$ are in $\zad{n}_{++}. \,$ 
We allow some or all components of $u$ to be $+ \infty$. If $u_i = + \infty,$ and $\alpha >0, \, $ then we define 
$\alpha u_i = + \infty, \,$ and if $b \in \zad{n}_{++}, \, $ is a row vector, then we define $bu = + \infty$. 
We will assume  $0 < \beta_1 \leq \beta_2 <  au. $
\eass
Recall the definition of a decomposable knapsack problem from Definition \ref{basicdef}. 
For the data vectors $p \,$ and $r \,$ from which we construct $a$ we will occasionally (but not always) assume
\bass \label{ass-pu}
$\, p \, \in \zad{n}_{++}, \, r \in \zad{n}, $ $p \,$ is not a multiple of $r, \,$ and 
\beq
r_1/p_1 \leq \dots \leq r_n/p_n.
\eeq
\co{\beq
\frac{r_1}{p_1} \leq \dots \leq \frac{r_n}{p_n}.
\eeq}
\eass

\nin{\bf Examples \ref{ex1} and \ref{ex2} continued} The problems \eref{ex1-prob} and 
\eref{ex2-prob} are DKPs with
\beqast
p & = & ( \,\,\, 1, \,\,\,\,1  ), \\
r & = & ( \,\,\, 1, -1  ), \\
u & = & ( \,\,\, 6, 6  ), \\
M & = & 20, \\
a & = & pM + r \, = \, ( 21, \, 19 ),
\eeqast
and 
\beqast
p & = & e, \\
r & = & 0, \\
u & = & e, \\
M & = & 2, \\
a & = & pM + r = 2e,
\eeqast
respectively.

\paragraph[wid]{Width and Integer width} 
\bdef 
Given a polyhedron $Q$, and an integral vector $c$, the width and the integer width of $Q$
in the direction of $c$ are
\beqast
\width(c, Q) & = & \max \, \{ \, cx \, | \, x \in Q \, \} - \min \, \{ \, cx \, | \, x \in Q \, \}, \\
\iwidth(c, Q) & = & \lfloor \max \, \{ \, cx \, | \, x \in Q \, \} \rfloor - \lceil \min \, \{ \, cx \, | \, x \in Q \, \} \rceil + 1. \\
\eeqast
If an integer programming  problem is labeled by ${\rm (P)}, \,$ and $c \,$ is an integral vector,
then with some abuse of notation  we  
denote by $\width(c, {\rm (P)})$ the width of the LP-relaxation of ${\rm (P)}$ in the 
direction $c$, and the meaning of $\iwidth(c, {\rm (P)})$ is similar.
\Edef
The quantity $\iwidth(c, Q) \,$ is the number of nodes generated by \bb \ when branching on the constraint $cx$.

\paragraph[br]{Basis Reduction} Recall the definition of a lattice generated by the columns of a rational
matrix $A$ from \eref{def-latt-B}. Suppose 
\beq
B = [ b_1, \dots, b_n ],
\eeq
with $b_i \in \zad{m}$. Due to the nature of our application, we will generally have 
$n \leq m. \,$  While most results in the literature are 
stated for full-dimensional lattices, it is easy to see that they actually apply to the general case.
Let $b_1^*, \dots, b_n^*$ be the Gram-Schmidt orthogonalization of $b_1, \dots, b_n, $ that is
\beq \label{bibist}
b_i =  \sum_{j=1}^{i} \mu_{ij} b_j^*,
\eeq
with 
\beq \label{muij} 
\ba{rcll} 
\mu_{ii} & = & 1 & \,\,\, (i=1, \dots, n), \\
\mu_{ij} & = & b_i^T b_j^*/\norm{b_j^*}^2 & \,\,\, (i=1, \dots, n; \, j = 1, \dots, i-1).
\ea
\eeq
We call $b_1, \dots, b_n$ {\em LLL-reduced} if 
\beqa \label{mucond}
| \mu_{ij} | & \leq & \frac{1}{2} \; (1 \leq j < i \leq n), \\ \label{exch-cond}
\norm{\mu_{i, i-1} b_{i-1}^* + b_i^*}^2 & \geq & \frac{3}{4} \norm{b_{i-1}^*}^2. 
\eeqa
An LLL-reduced basis can be computed in polynomial time for varying $n$.

Define the truncated sums
\beq
b_i(k) =  \sum_{j=k}^{i} \mu_{ij} b_j^* \,\, (1 \leq k \leq i \leq n),
\eeq
and for $i=1, \dots, n \,$ let $L_i$ be the lattice generated by 
$$
b_i(i), b_{i+1}(i), \dots, b_n(i).
$$
We  call $b_1, \dots, b_n$ {\em Korkhine-Zolotarev reduced} (KZ-reduced for short)
if $b_i(i)$ is the shortest lattice vector in $L_i \,$ for all $i$. Since $L_1 = L \,$ and $\, b_1(1)=b_1, \,$
in a KZ-reduced basis the first vector is the shortest vector of $L$. Computing the shortest vector in a 
lattice is expected to be hard, though it is not known to be NP-hard. It can be done 
in polynomial time when the dimension is fixed, and so can be computing a KZ-reduced basis.

\bdef
\label{def-redfactor}
Given a BR method (for instance LLL, or KZ), 
suppose there is a constant $c_n$ dependent only on $n$ with the
following property: for all full-dimensional lattices $\latt{A}$ in $\zad{n}, \,$ and for all reduced 
bases $\{ \, b_1, \dots, b_n \, \}$ of $\latt{A}$, 
\beqa \label{bidi}
\max  \{ \, \norm{b_1}, \dots, \norm{b_i} \, \} & \leq &
c_n \, \max  \{ \, \norm{d_1}, \dots, \norm{d_i} \, \}
\eeqa
for all $i \leq n, \,$ and any choice of linearly independent $d_1, \dots, d_i \, \in \latt{A}$.
We will then call $c_n$ the {\em reduction factor} of the BR method.
\end{Definition} 

The  reduction factors of LLL- and KZ-reduction are  $2^{(n-1)/2}$ (see \cite{LLL82}) and 
$\sqrt{n}$ (see \cite{Sh87}), respectively. For KZ-reduced bases, \cite{Sh87} gives 
a better bound, which depends on $i$, but for simplicity, we use $\sqrt{n}$. 

The $k$th successive minimum of the lattice $\latt{A}$ is
\beqast
\Lambda_k (\latt{B}) & = & \min \, \{ \, t \, | \, \exists \, k \, 
\text{linearly independent vectors in } \, \latt{A} \, 
\text{with norm at most} \, t \, \}.
\eeqast
So \eref{bidi} can be rephrased as 
\beqa \label{bidi-succ}
\max  \{ \, \norm{b_1}, \dots, \norm{b_i} \, \} & \leq &
c_n \Lambda_i(\latt{A}) \; \text{for all} \; i \leq n.
\eeqa

\paragraph[other]{Other notation} 
Given an integral matrix $C$ with independent rows, the
{\em null lattice}, or {\em kernel lattice} of $C$ is
\beqa
\nlatt{C} & = & \{ \, v \in \zad{n} \, | \, Cv=0 \, \}.
\eeqa

For vectors $f, \, p, \, $ and $u$ we write
\beq \label{maxdef}
\ba{rcl}
\max(f, p, \ell, u ) & = & \max \, \{ \, fx \, | \, px \leq \ell, \, 0 \leq x \leq u  \, \}, \\
\min(f, p, \ell, u ) & = & \min \, \{ \, fx \, | \, px \geq \ell, \, 0 \leq x \leq u \, \}.
\ea
\eeq
\co{
In our estimates we will use the following easy-to-prove lower bound on binomial coefficients:
\begin{Fact} \label{binfact} 
Suppose that $n, \,$ and $k \,$ are positive integers, $k \leq n.$ Then 
$$
\bpx 
n \\ k 
\epx \, \geq \, \left( \dfrac{n}{k} \right)^k.
$$
\end{Fact}
}

Theorems \ref{summ1} and \ref{summ2} below give a sample of our results from the following sections.
The overall results of the paper are more detailed,  but Theorems 1 and 2 are a convenient sample to first look at.

\bth \label{summ1}
Let $p \in \zad{n}_{++}, \, r \in \zad{n},  \, $  and $k$ and $M$ integers with 
\beq \label{kMbound} 
\ba{rcl} 
 0 & \leq & k \, < \, \sum_{i=1}^n p_i, \\
 M & >  & 2 \sqrt{n}(\norm{r}+1)^2 \norm{p}+1. 
\ea
\eeq
Then there are $\beta_1, \,$ and $\beta_2$ integers that satisfy 
\beq \label{kMsqrt} 
kM + \sqrt{n} \norm{r} < \beta_1 \leq \beta_2 < - \sqrt{n} \norm{r} + (k+1)M,
\eeq
and for all such $(\beta_1, \beta_2)$ 
the problem \eref{kp2} with $a = pM+r \,$ and $u = e$ has the following properties:
\benum
\item Its infeasibility is proven by $px \leq k \vee px \geq k+1$. 
\item Ordinary \bb \ needs at least $2^{\ell(p,k)}$ nodes to prove its infeasibility regardless of the order in which the branching variables 
are chosen. (Recall the definition of $\ell(p,k)$ from \eref{defratio}). 
\item The infeasibility of its rangespace reformulation computed with KZ-reduction is proven at the rootnode by branching on the last variable.
\eenum
\enth

\bth \label{summ2}
Let $p, \, $ and $r$ be integral vectors satisfying Assumption \ref{ass-pu}, $k \,$ and $M$ integers with 
\beq
\ba{rcl}
k & \geq & 0, \\
M & > & \max \, \{ \, k r_n/p_n - kr_1/p_1 - r_1/p_1 +1, \, 2 \sqrt{n} \norm{r}^2 \norm{p}^2 \}. 
\ea
\eeq
Then there exists a $\beta \,$ integer such that 
\beq \label{beta-pre} 
k (M + r_n/p_n)  \, < \, \beta \, < \, (k+1) (M + r_1/p_1),
\eeq
and for all such $\beta \,$ the problem \eref{kp2-eq} with $a = pM+r $ has the following properties:
\benum
\item Its infeasibility is proven by $px \leq k \vee px \geq k+1$. 
\item Ordinary \bb \ needs at least 
$$
\bpx 
\lfloor k/\norm{p}_\infty \rfloor + n-1  \\ n-1 
\epx \, 
$$
nodes to prove its infeasibility, independently of the sequence in which the 
branching variables are chosen. 
\item The infeasibility of the AHL reformulation computed with KZ-reduction is proven at the rootnode by branching on the last variable.
\eenum
\enth

\section{Why easiness for constraint branching implies  hardness for ordinary branch-and-bound}
\label{dkp-hard} 

In this section we prove a somewhat surprising result on instances of \eref{kp2}. If the infeasibility 
is proven by branching on $px, \,$ where $p$ is a positive integral vector, then this implies a lower bound on 
the number of nodes that ordinary \bb \ must take to prove infeasibility. 
So in a sense easiness implies hardness! 

A node of the branch-and-bound tree is identified by the subset of the variables that are fixed there, and by the values that they are fixed to.
We call $(\bx, F)$ a {\em node-fixing}, if  $F \subseteq N, \,$ and $\bx \in \zad{F} \,$ with 
$0 \leq \bx_i \leq u_i \, \forall i \in F, \,$ i.e. $\bx \,$ is a collection of integers corresponding to the components of $F$.

\bth \label{hardeasy} 
Let  $p \in \zad{n}_{++}, \,$  and $k$  an integer such  that the infeasibility of 
\eref{kp2} is proven by $\, px \leq k \vee px \geq k+1$. 
Recall the notation of $\ell(p,k)$ from \eref{defratio}. 

\benum

\item If $u = e, $ then ordinary \bb \ needs at least $2^{\ell(p,k)} \, $ nodes to prove the infeasibility of \eref{kp2}, independently of the sequence in which the 
branching variables are chosen. 

\item If $u_i = + \infty \, \forall i, $ then ordinary \bb \ needs at least 
$$
\bpx 
\lfloor k/\norm{p}_\infty \rfloor + n-1  \\ n-1 
\epx \, 
$$
nodes to prove the infeasibility of \eref{kp2}, independently of the sequence in which the 
branching variables are chosen. 
\eenum

\enth
\qed

To have a large lower bound on the number of \bb \ nodes that are necessary to prove
infeasibility, it is sufficient for  $\ell(p,k)$ to be large, which is true, if the components of $p$ are relatively small compared to $k, \,$ and are not too different. 
That is, we do not need the components of the {\em constraint vector} $a$ to be small, and not too different, as in Jeroslow's problem. 

First we need a lemma, for which one needs to recall the definition \eref{maxdef}.

\ble \label{technical} Let $k$ be an integer with $0 \leq k < pu. \,$ 
Then \eref{1} and \eref{2} below are equivalent: 
\benum
\item \label{1} The infeasibility of \eref{kp2} is proven by $\, px \leq k \vee px \geq k+1$.
\item \label{2} 
\beq \label{maxapkminapk} 
\max(a,p,k,u) < \beta_1 \leq \beta_2 < \min(a,p,k+1,u).
\eeq
\eenum
Furthermore,  if \eref{1} holds, then ordinary \bb \ cannot prune any node with node-fixing $(\bx,F)$ that satisfies 
\beq \label{partial} 
\sum_{i \in F} p_i \bx_i \leq k, \, \text{and} \, \sum_{i \not \in F} p_i u_i \geq k+1. 
\eeq
\ele
\pf{} Recall that we assume $0 < \beta_1 \leq \beta_2 < au. $ 
For brevity we will denote the box with upper bound $u$ by 
\beq \label{bu}
B_u = \{ \, x \, | \, 0 \leq x \leq u \, \}.
\eeq
The implication $\eref{2} \RA \eref{1}$ is trivial. To see $\eref{1} \RA \eref{2}$ first assume to the contrary that the lower inequality in 
\eref{maxapkminapk} is violated, i.e. there is $y_1 \,$ with 
\beq
y_1 \in B_u, \, p y_1 \leq k, \, a y_1 \geq \beta_1. 
\eeq
Let $x_1 = 0.$ Then clearly 
\beq
x_1 \in B_u, \, p x_1 \leq k, \, a x_1 < \beta_1. 
\eeq
So a convex combination of $x_1$ and $y_1, $ say $z$ satisfies 
\beq
z \in B_u, \, p z \leq k, \, a z  = \beta_1,
\eeq
a contradiction. Next, assume to the contrary that the upper inequality in 
\eref{maxapkminapk} is violated, i.e. there is $y_2 \,$ with 
\beq
y_2 \in B_u, \, p y_2 \geq k+1, \, a y_2 \leq \beta_2. 
\eeq
Define $x_2$ by setting its $i$th component to $u_i, \,$ if $u_i < + \infty, \,$ and to some large number $\alpha$ to be specified later, 
if $u_i = + \infty. \,$ 
If $\alpha$ is large enough, then 
\beq
x_2 \in B_u, \, p x_2 \geq k+1, \, a x_2 > \beta_2. 
\eeq
Then  a convex combination of $x_2$ and $y_2, $ say $w$ satisfies 
\beq
w \in B_u, \, p w \geq k+1, \, a w  = \beta_2,
\eeq
a contradiction. So $\eref{1} \RA \eref{2}$ is proven. 

Let $(\bx, F)$ be a node-fixing that satisfies \eref{partial}. Define $x'$ and $x''$ as 
\beq
x'_i \, = \, \left\{ \ba{rl} \bx_i & \text{if} \, i \in F \\
                            0     & \text{if} \, i \not \in F 
                    \ea
           \right.,  \,\, x''_i \, = \, \left\{ \ba{rl} \bx_i & \text{if} \, i \in F \\
                            u_i     & \text{if} \, i \not \in F. 
                    \ea
           \right. . \,\, 
\eeq
If $u_i = + \infty, \,$ then  $x''_i = u_i \,$ means ``set $x''_i$ to an $\alpha$ sufficiently large number''. 
We have $px' \leq k, \,$ so  $a x' < \beta_1;$ also, $px'' \geq k+1, \,$ so $ax'' > \beta_2$ holds as well.
Hence a convex combination of $x'$ and $x''$, say $z$ is LP-feasible for \eref{kp2}. Also, 
$z_i = \bx_i \,\, (i \in F) \,$ must hold, so the node with node-fixing $(\bx, F)$ is LP-feasible.

\qed

\pf{of Theorem \ref{hardeasy}}  
Again, we use the notation $B_u \,$ as in \eref{bu}. 

First we show that $0 \leq k < pu \,$ must hold. (The upper bound of course holds trivially, if any $u_i$ is $+ \infty$.) 
If $k < 0, \,$ then $px \geq k+1 \,$ is true for all $x \in B_u, \,$ so 
the infeasibility of \eref{kp2} could not be proven by $px \leq k \vee px \geq k+1.$ Similarly, if $k \geq pu, \,$ then 
$px \leq k \,$ is true for all $x \in B_u, \,$ so the infeasibility of \eref{kp2} could not be proven by $px \leq k \vee px \geq k+1.$ 

For both parts, assume w.l.o.g. that we branch on variables $x_1, x_2, \dots $ in this sequence. 
For part (1), let $F = \{ 1, 2, \dots, \ell(p,k) \}.$ From the definition of $\ell(p,k)$ it follows that any fixing of the variables in $F$ 
will satisfy \eref{partial}, so the corresponding node will be LP-feasible. Since there are $2^{\ell(p,k)}$ such nodes, the claim follows.

For part (2), let  $F = \{1, \dots, n-1 \},$ and assume that $x_i \,$ is fixed to $\bx_i \,$ for all $i \in F$. 
Since all $u_i$s are $+ \infty, \,$ 
this node-fixing will satisfy \eref{partial} if 
\beq \label{pibxi}
\bx_i \geq 0 \, \forall i \in F, \, \sum_{i \in F} p_i \bx_i \, \leq \, k.
\eeq
We will now give a lower bound on the number of $\bx \in \zad{F} \,$ that satisfy \eref{pibxi}.
Clearly, \eref{pibxi} holds, if 
\beq \label{pibxi-2}
\sum_{i \in F} \bx_i \, \leq \, \lfloor k/\norm{p}_{\infty} \rfloor
\eeq
does. It is known (see e.g., \cite{BR06}, page 30) that the number of nonnegative 
integral $(m_1, \dots, m_d) \,$ with  $m_1 + \dots + m_d \leq t$ is 
\beq \label{dt}
\bpx t + d  \\ d \epx. \, 
\eeq
Using this with $ t = \lfloor k/\norm{p}_{\infty} \rfloor, \, d = n-1, \,$ 
the number of $\bx \in \zad{F}\,$ that satisfy \eref{pibxi} is at least 
$$
\bpx 
\lfloor k/\norm{p}_\infty \rfloor + n-1  \\ n-1 
\epx,
$$
and so the number of LP feasible nodes is lower bounded by the same quantity.
\qed

\section{Recipes for decomposable knapsacks} 
\label{recipes}
\setcounter{equation}{0}

In this section we give simple recipes to find instances of \eref{kp2} and \eref{kp2-eq} with a decomposable structure. 
The input of the recipes is the $p$ and $r$ vectors, an integer $k, \,$ and the output is an integer $M, \, $ a vector $a$ with 
$a = pM+r, \,$ and the bounds $\beta_1$ and $\beta_2$, or $\beta$. 
The found instances will have their infeasibility proven by $px \leq k \vee px \geq k+1, \,$ and if $k$ is suitably chosen, 
be difficult for ordinary \bb \ by Theorem \ref{hardeasy}. We will show 
that several well-known hard integer programming instances are found by our recipes. 

\begin{figure}[h!]
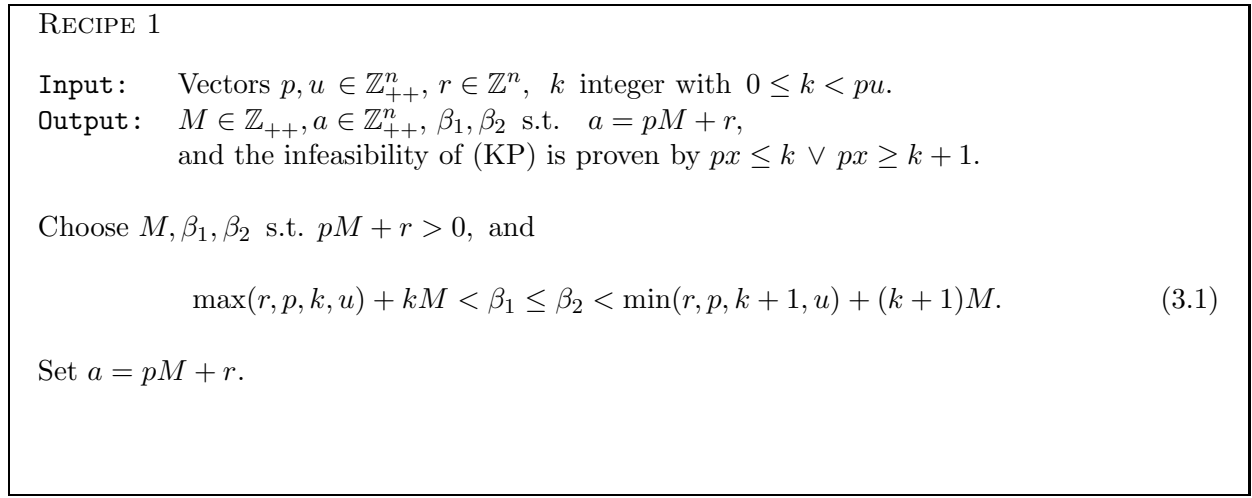

\framebox[6.5in]{\parbox{6.1in}{ 
{\sc Recipe 1}
\begin{tabbing}{}
*********\=******\=*******\=***\= \hspace{2.5in} \=  \kill  
{\tt Input:} \>   Vectors $p,  u \, \in \zad{n}_{++}, \, r \in \zad{n}, \, $  $\, k \,$ integer with $\, 0 \leq k < pu$. \\
{\tt Output:} \>   $M \in \zad{}_{++}, a \in \zad{n}_{++}, \, \beta_1, \beta_2 \,$ \st  $\,\,\, a = pM+r, \,$ \\
 \> and the infeasibility of \eref{kp2} is proven by $ px \leq k \, \vee \, px \geq k+1 $. \\
\\
Choose $M, \beta_1,  \beta_2 \,$ \st $pM + r > 0, \,$ and \\ \\
\> $ \,\,\, \max(r, p, k, u ) + k M < \beta_1 \leq \beta_2 < \min(r, p, k+1, u )  + (k+1)M. \hspace{.8in} \refstepcounter{equation} \label{maxrpkB} \, (\arabic{section}.\arabic{equation}) $ \\ \\
Set $a = pM+r$. \\ 
      \\
\end{tabbing}
}}
\caption{Recipe 1 to generate DKPs} 
\label{fig:dkp1}
\end{figure}

The recipes  are given in Figure \ref{fig:dkp1} and in Figure \ref{fig:dkp2}, respectively. 

\begin{figure}[h!]
\framebox[6.5in]{\parbox{6.1in}{ 
{\sc Recipe 2}
\begin{tabbing}{}
*********\=******\=*******\=***\= \hspace{2.5in} \=  \kill  
{\tt Input:} \>   Vectors $p, r \, \in \zad{n}\,$ satisfying Assumption  \ref{ass-pu}, $k$ nonnegative integer.  \\
{\tt Output:} \>   $M, \beta \in \zad{}_{++}, \, a \in \zad{n}_{++} \,$ \st $\,\, a = pM+r, \,$ and \\
 \> the infeasibility of \eref{kp2-eq} is proven by $px \leq k  \vee px \geq k+1. $ \\
\\
Choose $M, \beta \in \zad{}_{++} \,$ \st $pM+r > 0, \,$ and \\ \\
\> $ \di{\, 0 \leq  k (M + r_n/p_n)  \, < \, \beta \, < \, (k+1) (M + r_1/p_1).} \hspace{1.2in} \refstepcounter{equation} \label{beta-ineq} \, (\arabic{section}.\arabic{equation})  $
      \\ \\ 
Set $a = pM+r$. \\ 
\end{tabbing}
}}
\caption{Recipe 2 to generate instances of \eref{kp2-eq}} 
\label{fig:dkp2}
\end{figure}

\bth
Recipes 1 and 2 are correct. 
\enth
\pf{} 
Since $a = pM+r, \,$ 
\beqa \label{rec1}
\max(a,p,k,u) & \leq & \max(r,p,k,u) + kM,
\eeqa
and 
\beqa \label{rec2}
\min(a,p,k+1,u) & \geq & \min(r,p,k+1,u) + kM.
\eeqa
So the output of Recipe 1 satisfies 
\beq \label{maxapkminapk-2} 
\max(a,p,k,u) < \beta_1 \leq \beta_2 < \min(a,p,k+1,u),
\eeq
and so the infeasibility of the resulting DKP is proven by $px \leq k \vee px \geq k+1$. 

For Recipe 2, note that with the components of $u$ all equal to $+ \infty, \,$ we have 
\beqa
\max(r,p,k,u) & = & k r_n/p_n,  \, \text{and} \\
\min(r,p,k+1,u) & = & (k+1) r_1/p_1,  
\eeqa
so Recipe 2 is just a special case of Recipe 1. 
\qed

\nin{\bf Example \ref{ex1} continued}
We created Example \ref{ex1} using Recipe 1: 
here $pu = 12,$ so $k=5$ has $0 \leq k < pu$, and
\beqast
\max(r,p,k,u) & = & \max \, \{ \, x_1 - x_2 \, | \, 0 \leq x_1, x_2 \leq 6, \, x_1 + x_2 \leq 5 \, \} \, = \, 5, \,  \\
\min(r,p,k+1,u) & = & \min \, \{ \, x_1 - x_2 \, | \, 0 \leq x_1, x_2 \leq 6, \, x_1 + x_2 \geq 6 \, \} \, = \, -6.
\eeqast
So \eref{maxrpkB}   becomes
\beqast
5 + 5 M < \beta_1 \leq \beta_2 < -6 + 6M,
\eeqast
hence $M=20, \, \beta_1 = 106, \, \beta_2 = 113$ is a possible output of Recipe 1.

\nin{\bf Example \ref{ex2}  continued}
Example \ref{ex2} can also be constructed via Recipe 1:
now $pu =n$, so $k = (n-1)/2$ satisfies  $0 \leq k < pu$. Then $r=0$ implies
$$
\ba{rclcl}
\max(r,p,k,u) & = & \min(r,p,k+1,u) & = & 0,
\ea
$$
so
\eref{maxrpkB} becomes
\beqast
\frac{n-1}{2} M < \beta_1 \leq \beta_2 < \frac{n+1}{2} M,
\eeqast
and $M=2, \, \beta_1 = \beta_2  = n$ is a possible output of Recipe 1.
\qed

\bex \label{avistodd} 
{\rm 
Let $n$ be an odd integer, $k = \lfloor n/2 \rfloor, \, p  =   u  =  e, $ 
and $r$ an integral vector with 
\beqa
r_1 & \leq & r_2 \leq \dots \leq r_n.
\eeqa
Then we claim that any $M$ and $\beta = \beta_1 = \beta_2 $ is a possible output of Recipe 1, if 
\beq \label{betaM} 
\ba{rcl}
\beta & = & \, \lfloor (1/2) \sum_{i=1}^n (M+r_i) \rfloor, \\ 
M + r_1 & \geq & 0, \\ 
M + r_{k+1} & > & (r_{k+2} + \dots + r_n) - (r_1 + \dots + r_k).
\ea
\eeq
Indeed, this easily follows from 
\beqast
\max(r,p,k,u) & = &  r_{k+2} + \dots + r_n, \\
\min(r,p,k+1,u) & = & r_1 + \dots + r_k + r_{k+1}.
\eeqast
Two interesting, previously proposed  hard knapsack instances can be obtained by picking $r, \, M, \, $ and $\beta \,$ that satisfy 
\eref{betaM}. When 
\beq
\ba{rclrcl}
r & = & (2^{\ell+1}+1, \dots, 2^{\ell+n}+1), & M & = & 2^{n+\ell+1},
\ea
\eeq
with $\ell = \lfloor \log \, 2n \rfloor$, we obtain a feasibility version of a hard knapsack instance proposed by Todd in \cite{C80}. 
When 
\beq
\ba{rclrcl}
r & = & (1, \dots, n), & M & = & n(n+1),
\ea
\eeq
we obtain a feasibility version of a hard knapsack instance proposed by Avis in \cite{C80}. 

So the instances are 
$$
ax = \left\lfloor \dfrac{1}{2} \sum_{i=1}^n a_i \right\rfloor, \, x \in \{ 0, 1 \}^n, 
$$
with 
\beqa
a & = & (2^{n+\ell+1} + 2^{\ell+1}+1, \dots, 2^{n+\ell+1} + 2^{\ell+n}+1),
\eeqa
for the Todd-problem, and 
\beqa
a & = & ( n(n+1) + 1, \dots, n(n+1)+n )
\eeqa
for the Avis-problem.
}
\eex

\bex \label{reverse-avis} {\rm 
In this example we reverse the role of $p$ and $r$ from Example \ref{avistodd}, and will call the resulting DKP instance
a {\em reverse-Avis instance}.  This example illustrates how we can generate 
provably infeasible and provably hard instances from any $p$ and $r$;  also, the reverse-Avis instance 
will be harder from a practical viewpoint, as explained in Remark \ref{reverse} below. 

Let $n$ be a positive integer divisible by $4, \,$ 
\beq
\ba{rcl}
p & = & (1, \dots, n), \\
r & = & e, \\
k & = & n(n+1)/4.
\ea
\eeq
Since $k = (\sum_{i=1}^n p_i)/2, \,$ the first $3n/4$ components of $p$ sum to strictly more than $k, \,$ and the last $n/4$ sum to strictly less 
than $k, \,$ so 
\beq
\ba{rcl}
\max(r,p,k,u) & < & 3n/4, \\
\min(r,p,k+1,u) & > & n/4.
\ea
\eeq
Hence a straightforward computation shows that 
\beq
\ba{rcl}
M  & = & n/2 + 2, \\
\beta & = & \beta_1 = \beta_2 = 3n/4 + k(n/2+2) + 1
\ea
\eeq
are a possible output of Recipe 1. 
\co{
Also, $\ell(p,k) > n/4, \,$ so \eref{kp2} with $a = pM+r \,$ and the above $\beta$ has the following properties:
\benum
\item its infeasibility is proven by $px \leq k \vee px \geq k+1$; 
\item ordinary \bb \ needs at least $2^{n/4}$ nodes to prove the same.
\eenum
}
}
\eex

\bcor  \label{n2n4}
Ordinary \bb \ needs at least $2^{(n-1)/2}$ nodes to prove the infeasibility of Jeroslow's problem in Example \ref{ex2}, of the instances in Example 
\ref{avistodd} including the Avis- and Todd-problems, and at least $2^{n/4}$ nodes to prove the infeasibility of the reverse-Avis instance. 
\ecor
\pf{} We use Part (1)  of Theorem \ref{hardeasy}. In the first three instances $n$ is odd, $p = u = e, \, k = (n-1)/2, \, $ 
so $\ell(p,k) = k$. In the reverse-Avis instance we have $\ell(p,k) \geq n/4$ as explained after the definition \eref{defratio}. 

\qed

\brem \label{reverse} {\rm
While we can prove a $2^{(n-1)/2}$ 
lower bound for the Avis and the Todd instances, they are easy from a practical
viewpoint: it is straightforward to see that a single knapsack cover cut proves their infeasibility. 

For the reverse-Avis problem we can prove only a 
$2^{n/4}$ lower bound, but this problem is hard even from a practical viewpoint. We chose $n=60, $  and ran the resulting instance
using the CPLEX 11 MIP solver. After enumerating 10 million nodes the solver could not verify the infeasibility.
}
\erem

Next we give examples on the use of Recipe 2.
\bex \label{kp2-eq-ex}
{\rm 
Let $n=2, \,$ 
\beqast
k & = & 1, \\
p & = & ( \,\,\, 1, \,\,\,\,1  ), \\
r & = & ( \,\, -11, 5  ).
\eeqast
Then \eref{beta-ineq} in Recipe 2 becomes 
\beq
0 \leq M + 5 < \beta < 2(M-11),
\eeq
hence $M=29, \beta = 35 \,$ is a possible output of Recipe 2. So the infeasibility of 
\beq \label{18frob} 
\ba{rcl}
18 x_1 + 34  x_2 & = 35 \\
 x_1, x_2 & \in & \zad{}_+ 
\ea
\eeq
is proven by  $x_1 + x_2 \leq 1 \vee x_1 + x_2 \geq 2, \,$ a fact that is easy to check directly.
}
\eex
\qed

\bex \label{kp2-eq-ex-big}
{\rm In Recipe 2 $M$ and $\beta$ are constrained only by $r_1/p_1$ and $r_n/p_n$. 
So, if $n=17, \, k =1, \,$ and
\beqast
p & = & ( \,\,\, 1, 1, \dots, \,\,\,\,1  ), \\
r & = & ( \,\, -11, -10, \dots, 0, 1, \dots 5  ),
\eeqast
then $M=29, \beta = 35 \,$ is still a possible output of Recipe 2. So the infeasibility of 
\beq
\label{kp2-eq-ex-big-prob}
\ba{rcl}
18 x_1 + 19 x_2 + \dots + 34  x_2 & = 35 \\
 x_1, x_2, \dots, x_{17}  & \geq & 0 \\ 
 x_1, x_2, \dots, x_{17} & \in & \zad{}_+
\ea
\eeq
is proven by  $\sum_{i=1}^{17} x_i \leq 1 \vee \sum_{i=1}^{17} x_i \geq 2. \,$ 
}
\eex
\qed

We finally give an example, in which the problem data has polynomial size in $n, \,$ the infeasibility is proven by a split disjunction, but
the number of nodes that ordinary \bb \ must enumerate to do the same is a {\em superexponential} function of $n. \,$

\bex \label{kp2-eq-nt} {\rm
Let $n$ and $t$ be integers, $n, t \geq 2.$ We claim that the infeasibility of 
\beq
\label{nt-problem}
\ba{rcl}
(n^{t+1}+1) x_1 + \dots + (n^{t+1}+n)x_n & = & n^{2t+1} + n^{t+1} + 1, \\
                                    x_i  & \in & \zad{}_+ \, (i=1, \dots, n)
\ea
\eeq
is proven by 
\beq \label{ntsplit}
\sum_{i=1}^n x_i \leq n^t \, \vee \, \sum_{i=1}^n x_i \geq n^t+1,
\eeq
but ordinary \bb \ needs at least 
$$
n^{(n-1)(t-1)}
$$ 
nodes to prove the same. Indeed, 
\beq
\ba{rcl}
p & = & e, \\
r & = & (1, 2, \dots, n), \\
k & = & n^t, \\
M & = & n^{t+1}, \mbox{and}\\
\beta & = & n^{2t+1} + n^{t+1} + 1
\ea
\eeq
satisfy \eref{beta-ineq}. So the fact that the infeasibility is proven 
by \eref{ntsplit} follows from the correctness of Recipe 2. By Part (2) of Theorem \ref{hardeasy} ordinary \bb \ needs to enumerate at least
$$
\bpx 
n^t + n-1  \\ n-1 
\epx \, 
$$
nodes to prove the infeasibility of \eref{nt-problem}. But 
$$
\bpx 
 n^t + n-1  \\ n-1 
\epx \, \geq \, \bpx 
n^t  \\ n-1 
\epx \, \geq \, \bpx 
\dfrac{n^t}{n-1}
\epx^{n-1} \, \geq \, n^{(n-1)(t-1)}.
$$
}
\eex

\section{Large right hand sides in \eref{kp2-eq}. The branching Frobenius number} 
\label{sect-large-rhs}

In this section we assume that $p$ and $r$ integral vectors which satisfy Assumption \ref{ass-pu} are given, and 
let 
\beq
q = ( r_1/p_1, \dots, r_n/p_n).
\eeq
Recipe 2 returns a vector $a = pM+r, \,$ and an integral $\beta, \,$ such that the infeasibility of
\eref{kp2-eq} with this $\beta \,$ is proven by branching on $px$. 

The Frobenius number of $a$ is defined as the largest integer $\beta$ for which \eref{kp2-eq} is infeasible, and it is denoted 
by $\Frob(a). \,$ This section extends the lower bound result \eref{AL-lb} of Aardal and Lenstra in \cite{AL04, AL06} 
in two directions. First, using Recipe 2, we show that for sufficiently large $M$ there is a range of 
$\beta$ integers for which the infeasibility of 
\eref{kp2-eq} with $a = pM+r \,$ is proven by branching on $px$. The smallest such integer is essentially the same 
as the lower bound in  \eref{AL-lb}. 

We will denote 
\beq
f(M, \delta) = \left\lceil \frac{M + q_1 - \delta}{q_n - q_1} \right\rceil -1
\eeq
(for simplicity, the dependence on $p \,$ and $r \,$ is not shown in this definition).

\bth
\label{frob-thm} 
Suppose that $f(M,1) \geq 0 \, $ (i.e., $M \geq q_n - 2 q_1 + 1), \, a \in \zad{n}_{++}, \,$ with $a = pM+r.$ 
Then there is an integer $\beta \,$ with 
\beq \label{beta-ineq2}
f(M, 1) (M + q_n)  \, < \, \beta \, < \, (f(M, 1)+1) (M + q_1),
\eeq
and for all such $\beta$ integers the infeasibility of \eref{kp2-eq} is proven by 
$px \leq f(M, 1) \vee px \geq f(M,  1) +1 $. 
\enth
\pf{} There is an  integer $\beta$ satisfying \eref{beta-ineq} in Recipe 2, if 
\beq \label{kMq2}
k (M + q_n)  + 1  < (k +1) \left( M + q_1 \right).
\eeq
But it is straightforward to see that \eref{kMq2} is equivalent to $k \leq f(M,1). \,$ Choosing $k = f(M,1)$ turns \eref{beta-ineq} into 
\eref{beta-ineq2}. 
\qed

Clearly, for all $\beta \,$ right-hand sides found by Recipe 2
\beq \label{frob-bound-2}
\beta \leq \Frob(a).
\eeq
Since for the $\beta \,$ rhs values found by Recipe 2, the infeasibility of \eref{kp2-eq} has a short, split disjunction certificate, and 
there is no known ``easy'' method to 
prove the infeasibility of \eref{kp2-eq} with $\beta \,$ equal to $\Frob(a), \,$ 
such $\beta$ right-hand sides are interesting to study.

\begin{Definition}
Assume that $f(M,1) \geq 0, \,$ and $a$ is a positive integral vector of the form $a = pM+r. \,$ 
The {\em $p$-branching Frobenius number} of $a$ is the largest right-hand side for which 
the infeasibility of \eref{kp2-eq} is proven 
by branching on $px$. It is denoted by 
$$
\Frob_{p}(a).
$$
\end{Definition}

\bth
\label{frob-thm-2} 
Assume that $f(M,1) \geq 0, \,$ and $a$ is a positive integral vector of the form $a = pM+r. \,$ 
Then 
\beq \label{frobpx-bd} 
f(M, 1)(M + q_n)  < \Frob_{p}(a) < (f(M, 0)+1) \left( M + q_1 \right).
\eeq
\enth
\pf{} The lower bound comes from Theorem \ref{frob-thm}. 
Recall the notation \eref{maxdef}. If all components of $u$ are $ + \infty, \,$ then 
\beqa
\max(a, p, k, u) & = & k a_n/p_n \, = \, k(M+q_n),  \\
\min(a, p, k+1, u) & = & (k+1) a_1/p_1 \, = \, (k+1)(M+q_1).
\eeqa
So Lemma \ref{technical} implies that if the  infeasibility of \eref{kp2-eq} is proven by $px \leq k \vee px \geq k+1, \,$ then 
\beq \label{kMq1}
k (M + q_n)   <    \beta  < (k +1) \left( M + q_1 \right),
\eeq
hence 
\beq \label{kMq2-2}
k (M + q_n)   <    (k +1) \left( M + q_1 \right),
\eeq
which is equivalent to 
\beq \label{kMq3-2}
k   < \frac{M + q_1}{q_n - q_1} \, \LRA \, k   < \left\lceil \frac{M + q_1}{q_n - q_1} \right\rceil \, \LRA \, k \leq f(M,0).
\eeq
The infeasibility of \eref{kp2-eq} is proven by 
branching on $px \,$ iff it is proven by 
$px \leq k \, \vee \, px \geq k+1 \,$ for some nonnegative integer $k$. So, the largest such 
$\beta$ is strictly less than 
\beq \label{beta-ineq-3} 
(k +1) \left( M + q_1 \right),
\eeq
with $k  \leq f(M,0), \,$  so it is strictly less than 
$( f(M,0)+1) \left( M + q_1 \right), \,$ as required.

\qed

\nin {\bf Example \ref{kp2-eq-ex} continued} Recall that in this example 
\beqast
p & = & ( \,\,\, 1, \,\,\,\,1  ), \\
r & = & ( \,\, -11, 5  ),
\eeqast
so we have 
$q_1 = -11, \, q_2 = 5.$
So if $M=29, \,$ then $f(M,0) = f(M,1) = 1$, and the bounds in Theorem \ref{frob-thm} become 
$$
34 \, < \, \beta \, < \, 36.
$$
Hence Theorem \ref{frob-thm} finds only $\beta = 35, \,$ as the only integer for which the 
infeasibility of \eref{18frob} is proven by branching on $x_1 + x_2 \leq 1 \vee x_1 + x_2 \geq 2. \,$ 

Letting $a = pM + r = (18, 34), \,$ Theorem \ref{frob-thm-2} shows 
$$
34 \, < \, \Frob_p(a) \, < \, 36,
$$
so $\Frob_{p}(a) = 35. \, $ 
\co{
On the other hand, by a result of Brauer and Schockley, (\cite{BS62}), 
$
\Frob(a) \, = \, a_1 a_2 - a_1 - a_2 \, = \, 560,
$ so 
}

\section{The geometry of the original set, and the reformulation}
\label{geometry} 
\setcounter{equation}{0}

This section proves some basic results on the geometry of the reformulations using  
ideas from the recent article of Mehrotra and Li \cite{ML04}. Our goal is to relate the width
of a polyhedron to the width of its reformulation in a given direction.

\bth \label{range-geometry}
Let 
\beqast
Q & = &      \{ \, x \, \in \rad{n} \, | \, Ax \leq b \, \}, \\
\tQ & = &    \{ \, y \, \in \rad{n} \, | \, AUy \leq b \, \}, 
\eeqast
where $U$ is a unimodular matrix, and $c \in \zad{n}$. 

\nin Then 
\benum
\item \label{range-geometry-1} $$ \max \, \{ \, cx \, | \, x \in Q \, \}  \, = \, \max \, \{ \, cUy \, | \, y \in \tQ \, \}, $$
with $x^*$ attaining the maximum in $Q$ if and only if $U^{-1}x^*$ attains it in 
$\tQ$. 
\item \label{range-geometry-3} $$\width(c, Q) \, = \, \width(cU, \tQ). $$
\item \label{range-geometry-4} $$\iwidth(c, Q) \, = \, \iwidth(cU, \tQ). $$

\eenum

\enth
\pf{} Statement \eref{range-geometry-1} follows from 
\beq
Q \, = \,  \{ \, Uy \, | \, y \in \tQ \, \},
\eeq
and an analogous result holds for ``min''. Statements \eref{range-geometry-3}  and \eref{range-geometry-4} 
are easy consequences.
\qed

\nin Theorem \ref{range-geometry} immediately implies
\bcor \label{range-geometry-cor}
$$ \min_{{c \in \zad{n} \setminus \{ \, 0 \, \}}} \, \width(c, Q) \, = \, \min_{{d \in \zad{n}\setminus \{ \, 0 \, \} }} \, \width(d, \tQ). \,$$
\ecor
\qed

\bth \label{null-geometry} Suppose that the integral matrix $A$ has $n$ columns, and $m$ linearly independent rows, let
$S$ be a polyhedron, and 
\beqast
Q & = &      \{ \, x \in \rad{n} \, | \, x \in S, \, Ax = b \, \}, \\
\hQ & = &    \{ \, \lambda \, | \,  V \lambda + x_b \in S, \, \lambda \in \rad{n-m} \,  \,\}, 
\eeqast
where  $V$ is a basis matrix for $\nlatt{A}, \,$ and $x_b \,$ satisfies $A x_b = b. \, $ If $c  \, \in \zad{n}$ is a row vector, then 
\benum
\item \label{null-geometry-1} $$ \max \, \{ \, cx \, | \, x \in Q \, \}  \, = \, c x_b + \max \, \{ \, cV \lambda \, | \, \lambda \in \hQ \, \}, $$
with $x^*$ attaining the maximum in $Q$ if and only if $\lambda^*$ attains it in 
$\hQ$, where $x^* = V\lambda^* + x_b$.

\item \label{null-geometry-2} $$\width(c, Q) \, = \, \width(cV, \hQ). $$

\item \label{null-geometry-3} $$\iwidth(c, Q) \, = \, \iwidth(cV, \hQ). $$

\eenum

\enth
\pf{} Statement \eref{null-geometry-1} follows from 
$$
Q \, = \, \{ \, V \lambda + x_b \, | \, \lambda \in \hQ \}.
$$
An analogous result holds for ``min'', and  statements \eref{null-geometry-2} and \eref{null-geometry-3} are then straightforward consequences.

\qed

Theorem \ref{null-geometry} can be ``reversed''. That is, given a row vector $d \in \zad{n-m}$, we can find a row vector $c \in \zad{n}$, such that 
$$ 
\max \, \{ \, cx \, | \, x \in Q \, \}  \, = \, \max \, \{ \, d \lambda \, | \, \lambda \in \hQ \, \} + {\rm const}.
$$
Looking at  \eref{null-geometry-1} in Theorem \ref{null-geometry}, for the given $d$ it suffices to solve 
\beq \label{cvd}
cV = d, \, c \in \zad{n}.
\eeq
The latter task is trivial, if we have a $V^*$ integral matrix such that 
\beq \label{vstar}
V^*V = I_{n-m};
\eeq
then  $c = dV^*$ will solve \eref{cvd}. 
To find $V^*, \,$ let $W$ be an integral matrix such that $U = [W, V]$ is unimodular; for instance 
$W$ will do, if 
$$
A [W, V] = [H, 0],
$$
where $H$ is the Hermite Normal Form of $A$. Then we can choose $V^*$ as the submatrix of $U^{-1}$ consisting 
of the last $n-m$ rows.

In this way we have proved Theorem \ref{null-geometry-reverse}
and  Corollary \ref{null-geometry-reverse-cor}, which are 
essentially the same as Theorem 4.1, and Corollary 4.1 proven by Mehrotra and Li in 
\cite{ML04}:

\bth \label{null-geometry-reverse} (Mehrotra and Li) Let $Q, \, \hQ, \, V \, $ 
be as in Theorem  \ref{null-geometry}, and $V^*$ 
a matrix satisfying \eref{vstar}. Then 
\benum
\item \label{null-geometry-reverse-1} $$ \max \, \{ \, d \lambda \, | \, \lambda \in \hQ \, \} \, = \, \max \, \{ \, dV^* x \, | \, x \in Q \, \}  - dV^*x_b, $$
with $x^*$ attaining the maximum in $Q$ if and only if $\lambda^*$ attains it in 
$\hQ$, where $x^* = V\lambda^* + x_b$.

\item \label{null-geometry-reverse-2} $$\width(dV^*, Q) \, = \, \width(d, \hQ). $$

\item \label{null-geometry-reverse-3} $$\iwidth(dV^*, Q) \, = \, \iwidth(d, \hQ). $$

\eenum

\enth
\qed

\bcor (Mehrotra and Li) \label{null-geometry-reverse-cor} Let $Q, \hQ, V, \,$ and $V^*$ be as before. Then 
$$ \min_{{d \in \zad{n-m}\setminus \{ \, 0 \, \} }} \, \width(d, \hQ) \, = \, \min_{{c \in \latt{V^{*T}} \setminus \{ \, 0 \, \}}} \, \width(c, Q). \,$$
\ecor
\qed

\section{Why the reformulations make DKPs easy} 
\label{dkp-easy} 
\setcounter{equation}{0}

This section will assume a decomposable structure on \eref{kp2} and \eref{kp2-eq}, that is
\beqa
a & = & pM + r,
\eeqa
with $p \in \zad{n}_{++}, \, r \in \zad{n}, \,$ and $M$ an integer. We 
show that for large enough $M \,$ the phenomenon of Examples
\ref{ex1} and \ref{ex2} {\em must} happen, i.e., the originally difficult DKPs will turn 
into easy ones. 

We recall that for a given a matrix $A$, we use a Matlab-like notation, 
and denote its $j^{th}$ row, and column by 
$
A_{j, :} \; \text{and} \; A_{:, j},
$
respectively. 

An outline of the results is: 
\benum

\item If $M$ is large enough, and $U$ is the transformation matrix of the rangespace reformulation, then $pU$ will have 
a ``small'' number of nonzeros. Considering the equivalence between the old and new variables 
$Uy = x, \,$ this means that branching on just a few variables in the reformulation will ``simulate'' 
branching on the backbone constraint $px \,$ in the original problem. An analogous  result will hold 
for the AHL reformulation.

\item It is interesting to look at what happens, when branching on $px$ does not prove 
infeasibility in the original problem, but the width in the direction of $p$ is
relatively small --  
this is the case in \eref{kp2-eq} as we prove in Lemma \ref{smallwidth} below. 

Invoking the results in Section \ref{geometry} will prove that when
$M$ is sufficiently large, the same, or smaller width is achieved
along a unit direction in either one of the reformulations.

\eenum

\ble 
\label{smallwidth}
Suppose that Assumption \ref{ass-pu} holds. Then 
\beqa \label{width-p}
\width(p, {\mbox{\rm \eref{kp2-eq}}})  & = & \Theta(\beta/M^2), \\
\label{width-ei}
\width(e_i, {\mbox{\rm \eref{kp2-eq}}}) & = & \Theta(\beta/M) \; \forall i \in \{1, \dots, n \}.
\eeqa
In both equations  the constant  depends on $p \, $ and $r$.
\ele
\pf{}: Since $a_i = p_i M + r_i, \, $ 
\beq
r_1/p_1 \leq \dots \leq r_n/p_n
\eeq
implies 
\beq
p_1/a_1 \geq \dots \geq p_n/a_n.
\eeq
So
\beqast
\max \{ \, px \, | \, ax= \beta, \, x \geq 0 \, \} & = &  \beta p_1/a_1, \\
\min \{ \, px \, | \, ax= \beta, \, x \geq 0 \, \} & = & \beta p_n/a_n,
\eeqast
and therefore
\beqast
\width(p, {\mbox{\rm \eref{kp2-eq}}}) & = & \beta \, (  p_1/a_1 - p_n/a_n ) \\ 
                                                   & = & \beta ( p_1 a_n - p_n a_1)/(a_1 a_n) \\
                                                   & = & \beta ( p_1 r_n - p_n r_1)/(a_1 a_n). \\
\eeqast
Also, 
\beqast
\max \{ \, x_i \, | \, ax= \beta, \, x \geq 0 \, \} & = &  \beta/a_i, \\
\min \{ \, x_i \, | \, ax= \beta, \, x \geq 0 \, \} & = & 0,
\eeqast
hence 
\beqast
\width(e_i, {\mbox{\rm \eref{kp2-eq}}}) & = & \beta/a_i.
\eeqast
Since 
$$
a_i = \Theta(M) \, \, \forall i \in \{ \, 1, \dots, n \, \},
$$
both \eref{width-p} and \eref{width-ei} follow.
\qed

\subsection{Analysis of the rangespace reformulation}
\label{subsec-range}

\nin After the rangespace reformulation is applied, the problem \eref{kp2}
becomes
 \beq  \label{kptref} \tag{KP-R}
  \ba{rcl}
  \beta_1 \, \leq & (aU)y & \,\, \leq \beta_2 \\
  0 \, \leq & Uy & \,\, \leq u \\
   &  y & \, \in \zad{n},
  \ea
 \eeq
where the matrix $U$ was computed by a BR algorithm with input
\beq \label{defA}
A = \bpx
a \\
I
\epx = \bpx
pM + r \\
I
\epx.
\eeq
Let us write 
$$
\tA = AU, \, \ta = aU, \, \tp = pU, \, \tr = rU, 
$$ and fix
$c_n$, the reduction factor  of the used BR algorithm. 

Recall that for a lattice $L, \,$ 
$\Lambda_k(L) \,$ is the smallest real number $t$ for which there are $k \,$ linearly 
independent vectors in $L$ with norm at most $t$.

For brevity, we will denote 
\beqa \label{denote-alpha-k}
\alpha_k & = & \Lambda_k(\nlatt{p}) \; (k=1, \dots, n-1).
\eeqa

First we need a technical lemma:
\ble \label{succ-min-lemma-range}
Let  $A$ be as in \eref{defA}. 
Then 
\beqa \label{succ-range}
\Lambda_k ( \latt{A}) & \leq & ( \norm{r}+1) \alpha_k \; \text{for} \; k \in \{1, \dots, n-1 \}.
\eeqa
\ele
\pf{} 
We need to show that there are $k \,$ linearly independent vectors 
in $\latt{A} \,$ with norm bounded by $ (\norm{r}+1) \alpha_k. \,$

Suppose that $w_1, \dots, w_k \;$ are
linearly independent vectors in $\; \nlatt{p}$ with norm bounded by 
$\alpha_k$.   Then $A w_1, \dots, A
w_k \;$ are linearly independent in $\,\latt{A}$, and 
$$ A w_i \, = \, \bpx a \\ I \epx w_i \, = \, 
\bpx pM + r \\ I \epx w_i \,= \, \bpx
r w_i \\ w_i \epx \; \forall i, 
$$ 
hence
$$
\norm{ A w_i} \leq ( \norm{r}+1) \norm{w_i} \, \leq \, ( \norm{r}+1) \alpha_k \; (i=1, \dots, k)
$$
follows, which proves \eref{succ-range}.

\qed


\bth \label{range-main-thm} The following hold: 
\benum
\item \label{range-main-thm-1} Let $k \leq n-1, \,$ and suppose 
\beq \label{M-bound-k}
M > c_n (\norm{r}+1)^2 \alpha_k.
\eeq
Then 
 \beq \label{tp1k}
\tp_{1:k} = 0.
\eeq

Also, if the infeasibility of \eref{kp2} is proven by branching on $px, \,$ then the 
infeasibility of \eref{kptref} is proven by branching on $y_{k+1}, \dots, y_n$.

\item \label{range-main-thm-2} Suppose
\beq \label{M-bound-n}
M >  c_n (\norm{r}+1)^2 \norm{p}.
\eeq
Then 
\beq \label{tp1n}
\tp_{1:n-1} = 0,
\eeq
and
\beq \label{en-width}
\ba{rcl}
\width (e_n, {\mbox{\rm \eref{kptref}}})  & \leq & \width (p, {\mbox{\rm \eref{kp2}}}) \\
\iwidth (e_n, {\mbox{\rm \eref{kptref}}})  & \leq & \iwidth (p, {\mbox{\rm \eref{kp2}}}) 
\ea
\eeq

\nin In particular, in the rangespace reformulation of \eref{kp2-eq} the width, and the integer width in the direction of 
$e_n$ are
$$
\Theta(\beta/M^2).
$$

\eenum
\enth
\qed


Before proving Theorem \ref{range-main-thm}, we give some intuition to the validity of \eref{tp1k}, and \eref{tp1n}.
Suppose $M$ is ``large'', compared to $\norm{p}$, and $\norm{r}$.
In view of how the matrix $A$ looks in \eref{defA}, it is clear that its columns are 
{\em not} short, and near orthogonal, due to the presence of the nonzero $p_i$ components.
Thus to make its columns short and nearly orthogonal, the best thing to do is to apply 
a unimodular transformation that eliminates ``many'' nonzero $p_i$s. 

\nin\pf{} For brevity, denote by $Q\,$ and $\tQ\,$ the feasible set of the LP-relaxation
of \eref{kp2} and \eref{kptref}, respectively. 

\nin\pf{of \eref{range-main-thm-1}}  To show \eref{tp1k}, fix $j \leq k$; we will prove $\tp_j = 0$. 

\nin Since $\tA$ was computed by a BR algorithm with reduction factor $c_n$, Lemma \ref{succ-min-lemma-range}
implies 
 \beq \label{ta-fpr}
     \ba{rcl}
     \norm{\tA_{:,j}}  & \leq &  c_n (\norm{r}+1) \alpha_k.
     \ea
     \eeq
    To get a contradiction, suppose $\tp_j \neq 0$. Then, since $\tp_j$ is integral,
      \beq \label{ta-mfpr}
       \ba{rcl}
       \norm{\tA_{:,j}} & \geq & |\ta_{j}| \\
                   & =     & | \tp_j M + \tr_j | \\
                   & \geq  & M - | \tr_j|.
       \ea
       \eeq
       Hence 
      \beq \label{contr}
       \ba{rcl}
       M  & \leq & \norm{\tA_{:,j}} + | \tr_j| \\ 
          & \leq & \norm{\tA_{:,j}} + \norm{r} \norm{U_{:,j}}, \\
          & \leq & \norm{\tA_{:,j}} + \norm{r} \norm{\tA_{:,j}}, \\
                   & =  & (\norm{r}+1) \norm{\tA_{:,j}} \\
                   & \leq  & c_n (\norm{r}+1)^2 \alpha_k,
       \ea
       \eeq
         with the second inequality coming from Cauchy-Schwarz,
         the third from $U_{:,j}$ being a subvector of $\tA_{:,j}$, and the fourth from 
         \eref{ta-fpr}. Thus, we obtained a contradiction to the choice of $M$, which proves 
         $\tp_j = 0$. 

Suppose now that the infeasibility of \eref{kp2} is proven by branching on $px. \,$ 
We need to show:
\beq
y_i  \, \in \, \zad{} \, \forall \, i \in \{ k+1, \dots, n \} \, \RA \, y \not\in \tQ.
\eeq
Let $y \in \tQ$. Then 
$$
Uy \, \in \, Q \, \RA \, pUy \, \not \in \, \zad{} \, \RA \, \tp_{k+1} y_{k+1} +\dots + \tp_n y_n  \, \not \in \, \zad{} \RA y_i \, \not \in \, \zad{} \,\,\,\,\, \text{for some} \, i \in \{\, k+1, \dots, n \, \},
$$
as required. \co{
The  first implication is true, 
since the infeasibility of $Q \,$ (that is, of \eref{kp2}) is proven by branching on $px$.
The second follows from $\tp_{1:k} = 0$, and the third is trivial. }

\nin\pf{of \eref{range-main-thm-2}} The  statement \eref{tp1n} follows 
from \eref{tp1k}, and the obvious fact, that $\alpha_{n-1} \leq \norm{p}, \,$ 
since there  are $n-1$ linearly independent vectors in $\nlatt{p}$ with norm bounded by $\norm{p}$. 

To see \eref{en-width}, we claim 
\beqast
\width(e_n, \tQ) & \leq & \width( \tp_n e_n, \tQ) \\
                 & = & \width( pU, \tQ) \\
                 & = & \width( p, Q).
\eeqast
Indeed, the inequality follows from  $\tp_n$ being a nonzero integer.
The  first equality comes from  \eref{tp1n}, 
and the second from \eref{range-geometry-1} in Theorem \ref{range-geometry}. The inequalities hold, even if 
we replace ``$\width$'' by ``$\iwidth$'', so this proves the second inequality in \eref{en-width}.

The claim about the width in the direction of $e_n$ follows from 
\eref{en-width}, and Lemma \ref{smallwidth}. 

\qed
\subsection{Analysis of the AHL-reformulation}
\label{subsec-null}

The technique we use to analyse the AHL reformulation is similar, 
but the bound on $M, \,$ which is necessary for the dominant $p$ direction 
to turn into a unit direction is different.   If $\beta_1 = \beta_2 = \beta, $ then the AHL reformulation of \eref{kp2} is
 \beq  \label{kptref-null} \tag{KP-N}
  \ba{rcl}
   0 \, \leq & V \lambda + x_\beta & \,\, \leq u \\
   &  \lambda & \, \in \zad{n-1},
  \ea
 \eeq
where the matrix $V$ is a basis of $\nlatt{a}$
computed by a BR algorithm, and $ax_\beta = \beta$.

\nin Let us write $\hp = pV, \, \hr = rV$ and recall the notation for $\alpha_k \,$ from \eref{denote-alpha-k}.
Again we need a lemma. 
\ble \label{succ-min-lemma-null}
Let $k \in \{ \, 1, \dots, n-2 \, \}.$ Then 
\beqa \label{succ-null} 
\Lambda_k ( \nlatt{p} \cap \nlatt{r} ) & \leq & 2 \norm{r} \alpha_{k+1}^2.
\eeqa
\ele
\pf{} We need to show that there are $k \,$ linearly independent vectors 
in $ \nlatt{p} \cap \nlatt{r}\,$ 
with norm bounded by $ 2 \norm{r} \alpha_{k+1}^2. \,$

Suppose that $w_1, \dots, w_{k+1} \;$ are
linearly independent vectors in $\; \nlatt{p}$ with norm bounded by 
$\alpha_{k+1}$. Let $W = [ \, w_1, \dots, w_{k+1} \, ], \,$ and 
$$
d = rW \in \zad{k+1}.
$$
Suppose w.l.o.g. that for some $t \in \, \{ \, 1, \dots, k+1 \, \}$
$$
d_1 \neq 0, \, \dots, \, d_t \neq 0, \, d_{t+1} \, = \, \dots \, = \, d_{k+1} \, = \, 0.
$$
Then
$$
d_2 w_1 - d_1 w_2, \, d_3 w_1 - d_1 w_3, \, \dots, \, d_{t} w_{1} - d_{1} w_t
$$
are $t-1$ linearly independent vectors in $ \nlatt{p} \cap \nlatt{r} \,$  with norm bounded by 
\beqast
2 \norm{d}_\infty \alpha_{k+1} & = & 2 \di{ \left( \max_{ \, i=1, \dots, t \, } \, | r w_i | \, \right)}   \alpha_{k+1} \\
                         & \leq & 2 \norm{r} \di{  \left( \max_{ \, i=1, \dots, t \, } \, \norm{w_i} \,\right)} \alpha_{k+1}  \\ 
                         & \leq & 2 \norm{r} \alpha_{k+1}^2.
\eeqast
The $k+1-t$ vectors
$$
w_{t+1}, \dots, w_{k+1}
$$
are obviously in $ \nlatt{p} \cap \nlatt{r}, \,$  with their norm obeying the same bound, and 
the two groups together are linearly independent.

\qed

\bth \label{null-main-thm} 
Suppose that $p$ and $r$ are not parallel. Then the following hold: 
\benum
\item \label{null-main-thm-1} Let $k \leq n-2, \,$ and suppose 
\beq \label{M-bound-k-null}
M > 2 c_n \norm{r}^2 \alpha_{k+1}^2.
\eeq
Then 
\beq \label{hp1k}
\hp_{1:k} = 0.
\eeq
Also, if  
the infeasibility of \eref{kp2} is proven by branching on $px, \,$ then the 
infeasibility of \eref{kptref-null} is proven by branching on $\lambda_{k+1}, \dots, \lambda_{n-1}$. 

\item \label{null-main-thm-2} Suppose 
\beq \label{M-bound-n-null}
M > 2 c_n \norm{r}^2 \norm{p}^2.
\eeq
Then 
\beq \label{hp1n}
\hp_{1:n-2} = 0,
\eeq
and 
\beq  \label{wid-en}
\ba{rcl} 
\width (e_{n-1}, {\mbox{\rm \eref{kptref-null}}})  \leq \width (p, {\mbox{\rm \eref{kp2}}}) \\
\iwidth (e_{n-1}, {\mbox{\rm \eref{kptref-null}}})  \leq \iwidth (p, {\mbox{\rm \eref{kp2}}}).
\ea
\eeq
\nin In particular, in the AHL reformulation of \eref{kp2-eq} the width, and the integer width 
in the direction of $e_{n-1}$ are
$$
\Theta(\beta/M^2).
$$

\eenum
\enth

\pf{} First note that $pV \neq 0, \,$  since $aV =0, \, pV = 0 \,$ 
implies $rV = 0, \,$ hence $p$ and $r$ would be parallel.
Also, for brevity, denote by $Q\,$ and $\hQ\,$ the feasible set of the LP-relaxation
of \eref{kp2} and \eref{kptref-null}, respectively. 

\nin\pf{of \eref{null-main-thm-1}} To show \eref{hp1k}, 
fix $j \leq k$; we will prove $\hp_j = 0$.  Suppose to the contrary that 
$\hp_j \neq 0, \,$ then its absolute value is at least $1$. Hence 
      \beq \label{a-mgpr}
       \ba{rcl}
       0 = |a V_{:,j}| & =  & | \hp_j M + \hr_j | \\
                   & \geq  & M - |\hr_j |.
       \ea
       \eeq
Therefore
        \beqast
        M  & \leq  & |\hr_j | \\
           & =     &  |r V_{:,j}| \\
                 & \leq & \norm{r} \norm{V_{:,j}} \\
                 & \leq & 2 c_n \norm{r}^2 \alpha_{k+1}^2. 
        \eeqast
Here the second inequality comes from Cauchy-Schwarz. The third is true, since 
the columns of $V$ are a reduced basis of $\nlatt{a} \, \subseteq \nlatt{p} \cap \nlatt{r},$ 
and by using Lemma \ref{succ-min-lemma-null}. 

\nin Suppose now, that the infeasibility of \eref{kp2} is proven by branching on $px. \,$ 
We need to show:
\beq
\lambda_i  \, \in \, \zad{} \, \forall \, i \in \{ k+1, \dots, n-1 \} \, \RA \, \lambda \not\in \hQ.
\eeq
Let $\lambda \in \hQ$. Then 
$$
V \lambda + x_\beta \, \in \, Q \, \RA \, p (V \lambda + x_\beta) \, \not \in \, \zad{} \, \RA \, \hp_{k+1} \lambda_{k+1} +\dots + \hp_{n-1} \lambda_{n-1} + p x_\beta  \not \in  \zad{} \, \RA \, $$ \vspace{-.3in} $$  \lambda_i \not\in \zad{} \,\,\, \text{for some} \, i \in \{\, k+1, \dots, n-1 \, \},
$$
as required.

\nin\pf{of \eref{null-main-thm-2}} The statement \eref{hp1n} 
again follows from the fact that there  are $n-1$ linearly independent
vectors in $\nlatt{p}$ with norm bounded by $\norm{p}$. 

We will now prove \eref{wid-en}. 
Since $\hp_{n-1}$ is an integer, its absolute value is at least $1$. Hence
\beqast
\width(e_{n-1}, \hQ) & \leq & \width( \hp_{n-1} e_{n-1}, \hQ) \\
                 & = & \width( pV, \hQ) \\
                 & = & \width( p, Q),
\eeqast
with first equality true because of \eref{hp1n},
and the second one due to \eref{null-geometry-2} in Theorem \ref{null-geometry}. 
The proof of the integer width follows analogously. 

The claim about the width in the direction of $e_{n-1}$ follows from 
\eref{wid-en}, and Lemma \ref{smallwidth}. 

\qed

\subsection{Proof of Theorems \ref{summ1} and \ref{summ2}} 

\pf{of Theorem \ref{summ1}} 
Recipe 1 requires 
\beqa
\max(r, p, k, u ) + k M & < \beta_1 \leq \beta_2 & < \min(r, p, k+1, u )  + (k+1)M.
\eeqa
Since now $u=e, \,$ both $\max(r,p,k,u)$ and $\min(r,p,k+1,u)$ are bounded by $\norm{r}_1 \, \leq \, \sqrt{n} \norm{r}$ in absolute value.
\co{
\beq
\ba{rcl}
\max(r,p,k,u) & \leq & \norm{r}_1 \, \leq \, \sqrt{n} \norm{r}, \\
\min(r,p,k+1,u) & \geq & - \norm{r}_1 \, \geq \, - \sqrt{n} \norm{r}.
\ea
\eeq}
So if $\beta_1 \,$ and $\beta_2$ satisfy \eref{kMsqrt}, then they are a possible output of Recipe 1, so 
the infeasibility of the resulting DKP is proven by $px \leq k \vee px \geq k+1$. If 
\beq \label{bla} 
M > 2 \sqrt{n} \norm{r} +1,
\eeq
then there is room in \eref{kMsqrt} for $\beta_1$ and $\beta_2$ to be integers. 
Theorem \ref{hardeasy} implies the lower bound on the number of nodes that ordinary \bb \ must enumerate 
to prove infeasibility. 

On the other hand, \eref{range-main-thm-2} in Theorem \ref{range-main-thm} with $c_n = \sqrt{n}$ implies that if 
\beq \label{M-bound-n-2}
M >  \sqrt{n} (\norm{r}+1)^2 \norm{p},
\eeq
then the infeasibility of the rangespace reformulation is proven by branching on the last variable.

Finally, the bound on $M$ in \eref{kMbound} implies both \eref{bla} and \eref{M-bound-n-2}. 

\qed

\pf{of Theorem \ref{summ2}} 
From the lower bound on $M, \,$ there is a $\beta \,$ integer that satisfies \eref{beta-pre}, 
and the fact that the resulting instance's infeasibility is proven by 
$px \leq k \vee px \geq k+1 \,$ follows from the correctness of Recipe 2. The lower bound on the number of nodes that ordinary \bb \ must enumerate to prove infeasibility follows
from Theorem \ref{hardeasy}. 

The fact that the infeasibility of the AHL reformulation is proven by branching on the last variable follows from 
\eref{null-main-thm-2} in Theorem \ref{null-main-thm} with $c_n = \sqrt{n}$.
\qed

\section{A computational study} 
\label{comp}
\setcounter{equation}{0}

The theoretical part of the paper shows that 
\bit
\item DKPs with suitably chosen parameters are hard for ordinary \bb, and easy for branching
on $px$, just like Examples \ref{ex1}, \ref{ex2}, and
\item both the rangespace, and AHL reformulations make them easy.
The key point is that branching on the last few variables 
in the reformulation simulates the effect of branching on $px \,$ in the original problem. 
\eit
We now look at the question whether these results translate into 
practice. The papers  \cite{ABHLS00,AHL00,LW02,AL04} tested 
the AHL-reformulation on the following instances: 
\bit
\item In  \cite{AHL00}, equality constrained knapsacks arising from practical applications.
\item In  \cite{ABHLS00}, the marketshare problems  \cite{CD98}.
\item In  \cite{LW02}, an extension of the marketshare problems.
\item In \cite{AL04} the instances of \eref{kp2-eq}, with the rhs equal to $\Frob(a)$.
\eit
Our tested instances are bounded DKPs both with equality and inequality constraints, 
and instances of \eref{kp2-eq}. 

In summary, we found the following.

\benum
\item On {\em infeasible} problems, 
both reformulations are effective in reducing the 
solution time of proving infeasibility.

\item They are also effective on {\em feasible } problems. 

In feasible problems a solution may be found by accident, so it is not clear
how to theoretically quantify the effect of various branching strategies, 
or the reformulations on such instances.

\item They are also effective on optimization versions of DKPs. 

\item When $\beta_1 = \beta_2, \,$ i.e. both reformulations are applicable, there is
no significant difference in their performance. 

\eenum

The calculations are done on a Linux PC with a 3.2 GHz
CPU. The MIP solver was CPLEX 9.0. For feasibility versions of integer
programs, we used the sum of the variables as a dummy objective function.
The basis reduction computations called the Korkhine-Zolotarev (KZ) subroutines from the Number
Theory Library (NTL) version 5.4 (see \cite{NTL}). 

We let $n=50, \,$ and first generate 
$10$ vectors $p, r \in \zad{n} \,$ with the components of $p \,$ uniformly distributed in 
$[1,10]$ and the components of $r \,$ uniformly distributed in $[-10,10]$. We use these 
ten $\,p, r \,$ pairs for all families of our instances.

Recall the notation that for $k \in \zad{}, \, u \in \zad{n}_{++} \,$, 
\beqast
\max(r, p, k, u) & = & \max \, \{ \, rx \, | \, px \leq k, \, 0 \leq x \leq u  \, \}, \mbox{and}\\
\min(r, p, k+1, u) & = & \min \, \{ \, rx \, | \, px \geq k+1, \, 0 \leq x \leq u  \, \}. \\
\eeqast

\subsection{Bounded knapsack problems with $u=e$}

We used Recipe 1 to generate $10$ 
difficult DKPs, with bounds on the variables, as follows:

\nin For each $p, r$ we let 
$$
u=e, \, M=10000, \, k = n/2 = 25, \ a \, = \, pM + r,$$
and set 
\beqast
\beta_1 & = & \lceil \max(r, p, k, u ) + k M \rceil, \\
\beta_2  & = & \lfloor \min(r, p, k+1, u )  + (k+1)M \rfloor. \\
\eeqast
By the choice of the data $\beta_1 \leq \beta_2 \,$ holds in all cases.
We considered the following problems using these $a, \, u, \, \beta_1, \, \beta_2:$ 

\bit
\item
The basic infeasible knapsack problem:

 \beq  \label{kp2-feas-new} \tag{DKP-INFEAS}
  \ba{rcl}
\beta_1 \, \leq & ax & \,\, \leq \beta_2 \\
0 \, \leq & x & \,\, \leq u \\
   &  x & \, \in \zad{n}.
  \ea
  \eeq

\item The optimization version: 

 \beq  \label{kp2-opt-new} \tag{DKP-OPT}
  \ba{rcl} \max & ax \\
\st & ax & \,\, \leq \beta_2 \\
0 \, \leq & x & \,\, \leq u \\
   &  x & \, \in \zad{n}.
  \ea
  \eeq

We denote by $\beta_a$ the optimal value, and will use $\beta_a$ for creating further instances.

\item The feasibility problem, with the rhs 
equal to $\beta_a$:

 \beq  \label{kp2-feas-max} \tag{DKP-FEAS-MAX}
  \ba{rcl}
& ax & \,\, = \beta_a \\
0 \, \leq & x & \,\, \leq u \\
   &  x & \, \in \zad{n}.
  \ea
  \eeq

\item The feasibility problem, with the rhs 
set to make it infeasible: 

 \beq  \label{kp2-infeas} \tag{DKP-INFEAS-MIN}
  \ba{rcl}
& ax & \,\, = \beta_a + 1 \\
0 \, \leq & x & \,\, \leq u \\
   &  x & \, \in \zad{n}.
  \ea
  \eeq

\eit

On the last two families  both reformulations are applicable. 

The results are in Table \ref{n50k25u1M10000}. In the columns marked 'R', and 'N' we display the 
number of \bb \ nodes taken by CPLEX after rangespace and AHL-reformulation was applied, respectively.
In the columns marked 'ORIG' we show the number of \bb \ nodes taken by CPLEX 
on the original formulation. 

Since the LP subproblems of these instances are easy to solve, 
we feel that the number of \bb \ nodes is a better way of comparing the 
performance of the MIP solver with and without the reformulation. 

We also verified that providing $px$ as a branching direction in the
original formulation makes these problems easy.  We ran CPLEX on the
original instances, after adding a new variable $z, \,$ and the
equation $z = px, \,$ to the formulation. The results with this option
are essentially the same as the results in the 'R' and 'N' columns.

\subsection{Bounded knapsack problems with $u=10e$}

We repeated the above experiment with $u =10e, \,$ but all other settings the same. That is,
using the same ten $\,p,r \,$ pairs, we let
$$
u=10e, \, M=10000, \, k = n/2 = 25, \ a \, = \, pM + r,$$
and set 
\beqast
\beta_1 & = & \lceil \max(r, p, k, u ) + k M \rceil, \\
\beta_2  & = & \lfloor \min(r, p, k+1, u )  + (k+1)M \rfloor, \\
\eeqast
then solved the instances \eref{kp2-feas-new} and \eref{kp2-opt-new},
\eref{kp2-feas-max} and \eref{kp2-infeas} as before.  The results are
in Table \ref{n50k25u10M10000}. The original formulations turned out
to be more difficult now, whereas the reformulated problems were just
as easy as in the $u=e \,$ case.

\clearpage

\begin{landscape}

\begin{table}[ht!] 
\begin{tabular}{|r||c|c|c||c|c||c|c||c|c|c||c|c|c||} \hline 
\ \  & \multicolumn{3}{|c||}{RHS values} & \multicolumn{2}{|c||}{\eref{kp2-feas-new}} & \multicolumn{2}{|c||}{\eref{kp2-opt-new} } & \multicolumn{3}{|c||}{\eref{kp2-feas-max}} & \multicolumn{3}{|c||}{\eref{kp2-infeas}}\\ \hline 
Ins & $\beta_1 $ & $\beta_2 $ & $\beta_a$ & R  & ORIG & R & ORIG  & R & N & ORIG & R & N & ORIG \\ \hline 
1  & 250040 & 259972 & 250039  &        1 &    4330785  &       10 &  7360728  &       10 &        1 &  1219304  &        1 &        1 &  3181671  \\ \hline 
2  & 250044 & 259979 & 250043  &        1 &   2138598  &       10 &  2329217  &        1 &        1 &    24130  &        1 &        1 &  1880980  \\ \hline 
3  & 250069 & 259973 & 250068  &        1 &   12480272  &       20 & 14006843  &       10 &        3 &    13800  &        1 &        1 & 11993912 * \\ \hline 
4  & 250034 & 259961 & 250033  &        1 &    1454260  &       10 &  2800898  &        1 &        5 &   555144  &        1 &        1 &  2531222  \\ \hline 
5  & 250037 & 259975 & 250036  &        1 &    4811440  &       10 &  6715586  &        1 &       10 &   155670  &        1 &        1 &  4131652  \\ \hline 
6  & 250038 & 259981 & 250037  &        1 &    3239982  &       10 &  2659752  &       10 &       10 &   283776  &        1 &        1 &  3155522  \\ \hline 
7  & 250085 & 259948 & 250084  &        1 &   11579118  &       10 & 14598901  &       10 &        1 &   107170  &        1 &        1 & 10871441 * \\ \hline 
8  & 250052 & 259961 & 250051  &        1 &    8659516  &       10 & 15440957  &       10 &        1 &   486255  &        1 &        1 &  8097370  \\ \hline 
9  & 250045 & 259984 & 250044  &        1 &    6393700  &       20 & 12520666  &        1 &       10 &    82455  &        1 &        1 &  6346153  \\ \hline 
10  & 250061 & 259968 & 250060  &        1 &  12244168  &       10 & 14848327  &       10 &        1 &     3600  &        1 &        1 & 11929161 * \\ \hline 
\end{tabular} \\ 
\vspace*{-0.2in}
\parbox{6in}{\caption{\label{n50k25u1M10000} DKPs with $n=50, k=25, u=e, M=10000$. '*': $1$ hour time limit exceeded}}
\vspace*{0.15in}
\end{table}


\begin{table}[hb!] 
\begin{tabular}{|r||c|c|c||c|c||c|c||c|c|c||c|c|c||} \hline 
\ \ \  & \multicolumn{3}{|c||}{RHS values} & \multicolumn{2}{|c||}{\eref{kp2-feas-new}}  & \multicolumn{2}{|c||}{\eref{kp2-opt-new} } & \multicolumn{3}{|c||}{\eref{kp2-feas-max}} & \multicolumn{3}{|c||}{\eref{kp2-infeas}} \\ \hline 
Ins & $\beta'$ & $\beta$ & $\beta_a$ & R &  ORIG & R & ORIG  & R & N & ORIG & R & N & ORIG \\ \hline 
1  & 250083 & 259719 & 250082  &              1   & 13204411   &        1 & 12927001   &        1 &        1 &  2571521   &        1 &        1 & 11968829 * \\ \hline 
2  & 250111 & 259779 & 250110  &             1   & 13674751   &        1 & 13369911   &        1 &        1 & 12441612 * &        1 &        1 & 11968829 * \\ \hline 
3  & 250156 & 259729 & 250155  &             1   & 10939735   &        1 & 13737652   &        1 &        1 &  1702224   &        1 &        1 & 10342918 * \\ \hline 
4  & 250098 & 259619 & 250097  &               1   & 14678404   &        1 & 12762803   &        1 &        1 &    25917   &        1 &        1 & 13480436 * \\ \hline 
5  & 250059 & 259759 & 250058  &               1   & 14128736   &        1 & 13464255   &        1 &        1 &  5829029   &        1 &        1 & 13070602 * \\ \hline 
6  & 250051 & 259799 & 250050  &               1   & 13979145   &       10 & 12310057   &        1 &        1 &   597113   &        1 &        1 & 13211779 * \\ \hline 
7  & 250206 & 259489 & 250205  &               1   &  8895772   &       10 &  8725886   &        1 &       10 & 10046297 * &        1 &        1 & 13211779 * \\ \hline 
8  & 250111 & 259619 & 250110  &               1   & 13198252   &        1 & 13799370   &        1 &        1 & 12235292 * &        1 &        1 & 13211779 * \\ \hline 
9  & 250081 & 259849 & 250080  &                1   & 13136603   &       10 & 13082057   &        1 &        1 &    18687   &        1 &        1 & 12448850 * \\ \hline 
10  & 250206 & 259689 & 250205  &               1   &  9251523   &       10 & 12947576   &        1 &        1 &  9692170 * &        1 &        1 & 12448850 * \\ \hline 
\end{tabular} \\ 
\vspace*{-0.13in}
\parbox{6in}{\caption{\label{n50k25u10M10000} DKPs with n=50, k=25, u=10, M=10000. '*': $1$ hour time limit exceeded}}
\end{table} 

\clearpage

\end{landscape}

\subsection{Equality constrained, unbounded knapsack problems}

In this section we consider instances of the type 

\beq  \label{kp2-eq-new} \tag\mbox{{KP-EQ}}
  \ba{rl}
ax & \,= \beta \\
x  & \, \geq 0 \\
x & \, \in \zad{n},
  \ea
  \eeq

We recall the following facts:

\bit
\item 
If we choose $M$ sufficiently large, and a $\beta$ integer satisfying
\beq \label{beta-ineq-new}
0 \leq  \left( \left\lceil \frac{M + q_1 - 1}{q_n - q_1} \right\rceil -1 \right)(M + q_n)  \, < \, \beta \, < \, \left\lceil \frac{M + q_1 - 1}{q_n - q_1} \right\rceil (M + q_1),  
\eeq
then the infeasibility of \eref{kp2-eq-new} is proven by branching on $px$.

\item If $\beta^*$ is the largest integer satisfying \eref{beta-ineq-new}, and 
$\Frob(a)$ the Frobenius number of $a, \,$ then clearly
$
\beta^* \leq \Frob(a).
$

\item Finding $\beta^*$ is trivial, while computing $\Frob(a) \,$ requires solving a sequence of integer programs.

\eit

We generated $20$ instances as follows: using the same $p, r \,$ pairs as 
in the previous experiments, we let 
$$
M = 10000.
$$
Then the first instance with a fixed $p, r \,$ pair arises by letting
the rhs in \eref{kp2-eq-new} be $\beta^*, \,$ and the second by
letting it to be equal to $\Frob(a)$.

The \eref{kp2-eq-new} instances with $\beta = \Frob(a)$ were already considered in \cite{AL04}.

Our computational results are in Table \ref{n50M10000}. Where we go
further than \cite{AL04} is by showing the following.
\bit

\item  $\beta^*$ and $\Frob(a)$ are not too different, and neither is the difficulty of \eref{kp2-eq-new} with 
these two different rhs values. 

\item Now both reformulations can be applied, and their performance is similar. 

\item According to Lemma \ref{smallwidth},  $\width(p, {\mbox{\rm \eref{kp2-eq}}}) \,$ and 
$\iwidth(p, {\mbox{\rm \eref{kp2-eq}}}) \,$ both should be small compared to 
the width in unit directions, even when the infeasibility of \eref{kp2-eq} is not proven 
by branching on $px$. This is indeed the case when $\beta = \Frob(a), \,$  
and we list $\iwidth(p, {\mbox{\rm \eref{kp2-eq}}}) \,$ in Table \ref{n50M10000} as well.

\item In the column ``$px$'' we list the number of \bb \ nodes
necessary to solve the problems, when the variable $z, \,$ and the
equation $z = px$, is added to the original problems.  The results are
similar to the ones obtained with the reformulations.

\eit


\begin{table}[ht] 
\begin{tabular}{|r||c|c|c||c|c|c|c||c|c|c|c||} \hline 
\ \ \  & \multicolumn{3}{|c||}{RHS, width} & \multicolumn{4}{|c||}{$ax=\beta^*$} & \multicolumn{4}{|c||}{$ax=\Frob(a)$ }\\ \hline 
\# & $\beta^*$ & ${\Frob(a)}$ & $\iwidth(p)$ & R & N & $px$ & ORIG & R & N & $px$ & ORIG \\ \hline 
1  & 7683078 & 7703088 &   1  &        1 &        1 &        1   &  7110020 \dag &        5 &        1 &       13   &  7320060 \dag \\ \hline 
2  & 8683916 & 8703917 &   1  &        1 &        1 &        1   &  6997704 \dag &        8 &        3 &       15   &  7123300 \dag \\ \hline 
3  & 8325834 & 8345840 &   1  &        1 &        1 &        1   & 15383299 \dag &        8 &        1 &       19   & 15313074 \dag \\ \hline 
4  & 10347239 & 10367238 &  2  &        1 &        1 &        1   & 10053497 \dag &       14 &       24 &       18   &  9928134 \dag \\ \hline 
5  & 16655001 & 16665004 &  1  &        1 &        1 &        1   &  9254836 \dag &       33 &       19 &       57   &  7519023 \dag \\ \hline 
6  & 9081818 & 9121828 &   1  &        1 &        1 &        1   &  6802797 \dag &       21 &        1 &       45   &  7011946 \dag \\ \hline 
7  & 6245624 & 6245632 &   1  &        1 &        1 &        1   &  7180978 \dag &        1 &        1 &       67   &  7151382 \dag \\ \hline 
8  & 10514739 & 10534740 &  1  &        1 &        1 &        1   &  7164967 \dag &        1 &        1 &       20   &  7178052 \dag \\ \hline 
9  & 14275715 & 14285716 &  1  &        1 &        1 &        1   &  7319379 \dag &        1 &        1 &       11   &  7368436 \dag \\ \hline 
10  & 9838851 & 9838851 &   0  &        1 &        1 &        1   &  7520143 \dag &        1 &        1 &        1   &  7230420 \dag \\ \hline 
\end{tabular} \\ 
\caption{\label{n50M10000} n=50, M=10000. '\dag': $1$ hour time limit exceeded}
\end{table} 

\subsection{Reformulated problems with basic MIP settings}

To confirm the easiness of the reformulated instances we reran all of
them with the most basic CPLEX settings: no cuts, no aggregator, no
presolve, and node selection set to depth first search. All instances
finished within a hundred nodes.

The instances and parameter files are publicly available from \cite{DKP}.

\section{Comments on the analysis in \cite{AL04}} 
\label{critique}
\setcounter{equation}{0}

In \cite{AL04}, Aardal and Lenstra studied the instances 
\eref{kp2-eq} 
with the constraint vector $a$ decomposing as 
\beq
a = pM+r, \, 
\eeq
with $p \in \zad{n}_{++}, r \in \zad{n}, \, M$ a positive integer, 
under Assumption \ref{al-assumption}. 
Recall that the reformulation \eref{eq-ip-ref} is constructed so that the columns of $B$ form an LLL-reduced basis of $\nlatt{a}$.

Denoting the last column of $B$ by $b_{n-1},\,$  Theorem 4 in \cite{AL04} proves \eref{AL-lb2}, which we recall here:
\beq \label{AL-lb2-new}
\norm{b_{n-1}} \geq \frac{\norm{a}}{\sqrt{ \norm{p}^2\norm{r}^2 - (pr^T)^2 }},
\eeq
and the following claims are made: 

\bit
\item
It can be assumed without loss of generality, that the columns of $B$ are ordered 
in a way that the first $n-2 \,$ form a basis for $\nlatt{p} \cap \nlatt{r}.$ This claim is used in the proof of Theorem 4.
\item Denoting by $Q$ the feasible set of the LP-relaxation of 
\eref{eq-ip-ref}, $b_{n-1}$ being long implies that $\iwidth(e_{n-1}, Q)$ is small. 
\eit
To reconcile the notation with that of \cite{AL04}, we remark that in the latter $L_0 \,$ and $L_C \,$ are used, where 
\beqast
L_0 & = & \nlatt{a}, \mbox{and}\\
L_C & = & \nlatt{p} \cap \nlatt{r}.
\eeqast

Here we provide Example \ref{alcounterex1}, in which $p, \, r, \, M \,$ satisfy 
Assumption \ref{al-assumption}, $B$ is an LLL-reduced basis of $\nlatt{a}, \,$ but $pB$ has $2$ nonzero components, so the first claim
does not hold. In Example \ref{alcounterex2}, using a modification of a construction of Kannan in \cite{K87agn}
we show a bounded polyhedron where the columns of the  constraint matrix are LLL-reduced, but branching on a variable corresponding to the longest column
produces exponentially many nodes. (Note that the polyhedron in \cite{AL04} is unbounded.) 
Finally, in Remark \ref{al-remark} we clarify the connection with our results. 

\bex \label{alcounterex1} {\rm 
Let $n=6, \, $ 
\beqast
p & = & (~1,~~~1,~~~3,~~~3,~~~3,~~~3), \\
r & = & (-7,-4,-11,-6,-5,-1), \\
M & = & 24, \\
a & = & (17,~~20,~~61,~66,~67,~71), \\
B & = & \bbx 
     ~1 & ~0 & -3 & ~1 & ~0 \\
     ~2 & -1 & -1 & -1 & ~0 \\
     -1 & -2 & ~0 & ~0 & -1 \\
     ~0 & ~0 & ~0 & ~1 & ~2 \\
     -1 & ~0 & ~0 & -2 & ~0 \\
     ~1 & ~2 & ~1 & ~1 & -1  \ebx, \\
\eeqast
The columns of $B$ form an LLL-reduced basis of $\nlatt{a}. \,$ They form a basis, since with 
$$
v = (0, -3, 1, 0, 0, 0)^T,
$$
the matrix $[B, v]$ is unimodular, and 
$$
a [ B, v] = [ 0_{1 \ti (n-1)}, \, \gcd(a) ].
$$
LLL-reducedness is straightforward to check using the definition. 
But 
\beqast
pB & = & (~0,~-1,~-1,~~0,~~0~),
\eeqast
so we cannot choose  $n-2=4$ columns of $B$ which would form a basis of $\nlatt{p} \cap \nlatt{r}.$
}
\eex

\bex \label{alcounterex2} {\rm 
Let $\rho$ be a real number in $(\sqrt{3}/2,1), \,$ and define the columns of the matrix $B \in \rad{n \ti n}$ as 
\begin{equation} 
\begin{array}{rcl}
b_1 & = & (\rho^0, 0, \dots, 0 )^T \\
b_2 & = & (\rho^0/2, \rho^1, \dots, 0, 0 )^T \\ 
b_3 & = & (\rho^0/2, \rho^1/2, \rho^2, \dots, 0, 0 )^T \\
    & \vdots &  \\
b_n & = & (\rho^0/2, \rho^1/2, \rho^2/2, \dots, \rho^{n-2}/2, \rho^{n-1})^T.
\end{array}
\end{equation}
Consider the polyhedron 
$$
Q \, = \, \{ \, \lambda \, | \, 0 \leq B \lambda \leq e_{n} \, \}.
$$
\bprop
The following hold:
\benum
\item \label{111} The columns of $B$ are an LLL-reduced basis of the lattice that they generate.
\item \label{222} $b_n$ is the longest among the $b_i$.
\item \label{333} 
$\width(e_{n}, Q)  > c^n \norm{b_n}$ for some $c > 1$.
\eenum
\eprop
\pf{} We have 
$$
b_i^* \, = \, \rho^{i-1}e_i \,\, (i=1, \dots, n),
$$
and when writing $b_i =  \sum_{j=1}^{i} \mu_{ij} b_j^*,$ 
\beq \label{mui2}
\mu_{i, i-1} = 1/2.
\eeq
Thus \eref{exch-cond} in the definition of LLL-reducedness becomes
\beq
\norm{b_i^*} \geq \frac{1}{\sqrt{2}} \norm{b_{i-1}^*},
\eeq
which follows from $\rho \geq 1/\sqrt{2}.$
Since 
\beq
\norm{b_n}^2 \, = \, \norm{b_{n-1}}^2 + \rho^{2(n-2)} \left( \rho^2 - \frac{3}{4} \right),
\eeq
this implies \eref{222}. 
By the definition of Gram-Schmidt orthogonalization, for any $\lambda_n \,$ we can always set the other $\lambda_i$ so 
$$
B \lambda = \lambda_n b_n^* = \lambda_n \rho^{n-1} e_{n},
$$
so 
$$
\width(e_{n}, Q) \geq \frac{2}{\rho^{n-1}}.
$$
and \eref{333} follows, since $\rho < 1, \, $ and $\norm{b_n} \leq n+1$.

We can make $B$ integral and still have \eref{111} through \eref{333} hold, by scaling, and rounding it. 
\qed

}
\eex

\brem{\rm 
\label{al-remark} 
\co{Under closer scrutiny, Theorem 4 in \cite{Al04} actually proves \eref{AL-lb2-new} with 
$\norm{b_{n-1}^*}$ instead of $\norm{b_{n-1}}$. }

Our Theorem \ref{null-main-thm} proves that if $M >  2^{(n+1)/2} \norm{r}^2 \norm{p}^2, \,$ and the reformulation is computed using 
LLL-reduction, then $(pB)_{1:(n-2)} = 0, \,$ and from this it does follow that the first $n-2$ columns 
of $B$ form a  basis of $\nlatt{p} \cap \nlatt{r}.$ 

Theorem \ref{null-main-thm} then finishes the analysis in a different way, 
by directly proving small width, namely showing 
\beq  \label{wid-en2}
\ba{rcl} 
\width (e_{n-1}, {\mbox{\rm \eref{kptref-null}}})  \leq \width (p, {\mbox{\rm \eref{kp2}}}), \mbox{and} \\
\iwidth (e_{n-1}, {\mbox{\rm \eref{kptref-null}}})  \leq \iwidth (p, {\mbox{\rm \eref{kp2}}}).
\ea
\eeq

}

\erem

\nin{\bf Acknowledgement} We thank the two referees, and the Associate Editor for their suggestions, which greatly improved the paper. 
Thanks are due to Mustafa Tural for his careful reading, and suggestions, and finding reference \cite{Brad71}.

\appendix
\label{app-sect}

\bibliography{IP_Refs}

\end{document}